\documentclass[11pt,letterpaper,reqno]{amsart}
\usepackage[T1]{fontenc}
\usepackage[utf8]{inputenc}
\usepackage{lmodern,amsmath,amssymb,amsthm,mathtools}
\usepackage[margin=1.25in]{geometry}
\usepackage{microtype,enumitem}
\usepackage{tikz,float}
\usetikzlibrary{arrows.meta,calc}
\usepackage[colorlinks=true,linkcolor=blue!50!black,urlcolor=blue!60!black,citecolor=blue!50!black]{hyperref}
\usepackage{xcolor}
\usepackage{xurl}
\usepackage[backend=bibtex,style=numeric,sorting=nyt,maxbibnames=99]{biblatex}
\AtEveryBibitem{%
  \iffieldundef{doi}{}{\clearfield{url}}%
  \iffieldequalstr{eprinttype}{arxiv}{\clearfield{url}}{}%
}
\newtheorem{theorem}{Theorem}[section]
\newtheorem{lemma}[theorem]{Lemma}
\newtheorem{proposition}[theorem]{Proposition}
\newtheorem{corollary}[theorem]{Corollary}
\theoremstyle{definition}
\newtheorem{remark}[theorem]{Remark}
\newcommand{\R}{\mathbb R}
\newcommand{\C}{\mathbb C}
\newcommand{\N}{\mathbb N}
\newcommand{\PP}{\mathbb P}
\newcommand{\EE}{\mathbb E}
\newcommand{\HS}{\mathrm{HS}}
\newcommand{\Tr}{\operatorname{Tr}}
\newcommand{\tr}{\operatorname{tr}}
\newcommand{\smin}{s_{\min}}
\newcommand{\dist}{\operatorname{dist}}
\newcommand{\rank}{\operatorname{rank}}
\newcommand{\diag}{\operatorname{diag}}
\newcommand{\alg}{\operatorname{alg}}
\newcommand{\Gr}{\operatorname{Gr}}
\newcommand{\one}{\mathbf 1}
\newcommand{\dd}{\,\mathrm d}
\newcommand{\calB}{\mathcal B}
\newcommand{\calD}{\mathcal D}
\newcommand{\frakD}{\mathfrak D}
\newcommand{\norm}[1]{\left\lVert#1\right\rVert}
\numberwithin{equation}{section}
\hypersetup{
 pdftitle={Universality of Brown measure convergence for non-Gaussian matrix polynomials with adjoints},
 pdfauthor={Yi Han},
 pdfsubject={Random matrices and free probability},
 pdfkeywords={Brown measure, Wigner matrices, least singular values, noncommutative polynomials}
}

\title[Universality of Brown measure convergence]{Universality of Brown measure convergence\\
for non-Gaussian matrix polynomials with adjoints}
\author{Yi Han}
\thanks{Department of Mathematics, Massachusetts Institute of Technology,
Cambridge, MA.\newline Email: \texttt{hanyi16@mit.edu}.}
\date{November 14, 2026}
\begin{document}
\begin{abstract}
Let $P$ be a noncommutative polynomial of arbitrary degree. We prove
that the empirical spectral distribution of $P$ evaluated at independent
Wigner matrices converges weakly in probability to its Brown measure
limit. The same conclusion holds for arbitrary $*$-polynomials in
independent iid matrices and their adjoints. The entry laws have bounded
real or joint planar densities and a finite $(2+\varepsilon)$ moment;
only second moments are required on the Wigner diagonal. The real and
imaginary parts of an entry may be dependent, and its pseudovariance
need not vanish. The proof develops the coefficient-algebra approach of
the author's earlier work on matrix polynomials with Gaussian and
non-Gaussian iid inputs \cite{Han}.
The main step is a polynomial lower-tail estimate, under only second
moments, for the least singular value of a tensor pencil with arbitrary
fixed matrix coefficients. We use a logarithmic determinant method:
a real-block small-ball argument shows that logarithmic derivatives of
the determinant control the inverse pencil. Relative polynomial
anti-concentration then gives the lower-tail bound. We establish
convergence of logarithmic potentials by Hermitian resolvent comparison
and rank truncation.
This extends the recent Brown measure convergence theorem of Collins, Kitagawa,
Miyagawa, Yuan and Zhong \cite{Collins} for GUE polynomials to non-Gaussian
Wigner inputs.
\end{abstract}
\maketitle
\section{Introduction}\label{sec:introduction}

For a matrix $M\in M_d(\C)$ with eigenvalues
$\lambda_1(M),\ldots,\lambda_d(M)$, counted with algebraic
multiplicity, its empirical spectral distribution is the probability
measure
\[
 \mu_M=\frac1d\sum_{j=1}^d\delta_{\lambda_j(M)}.
\]
A central problem in random matrix theory is to determine the limit of
$\mu_{M_N}$ for a sequence of random matrices $M_N$. For Hermitian
matrices this can often be done by computing traces of powers. For
nonnormal matrices, convergence of the $*$-moments alone does not
determine the limiting eigenvalue measure. Girko's Hermitization method
\cite{Girko} addresses this problem through the singular values of
$M_N-zI$. In particular, it requires control of their logarithms near
zero. The finite-variance circular law and its replacement principle
provide a basic example of this method \cite{TVK}; we refer to
\cite{BC} for a systematic account.

We study fixed noncommutative polynomials in independent iid matrices
and their adjoints, as well as polynomials in independent Wigner matrices.
For iid inputs, the entries of the polynomial need not be independent.
In both models, we prove convergence of the
empirical spectral distribution to the corresponding Brown measure
under bounded densities and a finite $(2+\varepsilon)$ moment.
The diagonal laws in the Wigner model need only a second moment.
We prove this using a least singular value estimate for tensor pencils
with arbitrary fixed coefficients. This estimate requires only second moments.

\subsection{Brown measure and notation}
For $M\in M_d(\C)$, set
\[
 U_M(z)=\frac1d\log|\det(M-zI_d)|.
\]
Then $\mu_M=(2\pi)^{-1}\Delta U_M$ in the sense of distributions.
The corresponding measure for an operator in a tracial von Neumann
algebra is the Brown measure \cite{Brown}, defined through the
Fuglede--Kadison determinant \cite{FK}. More precisely, let
$(\mathcal M,\tau)$ be a tracial von Neumann algebra with faithful
normalized trace. For $a\in\mathcal M$, write
\[
 U_a(z)=\tau(\log|a-z|),\qquad
 \mu_a=\frac1{2\pi}\Delta U_a.
\]
The logarithm is defined by decreasing positive regularization. Its
value at a particular point may initially be $-\infty$, whereas the
potential is locally integrable as a function of $z$.
We use standard free circular variables $c_1,\ldots,c_n$, normalized
by $\tau(c_\ell c_\ell^*)=1$, and standard free semicircular variables
$s_1,\ldots,s_n$, normalized by $\tau(s_\ell^2)=1$.
For free independence and circular and semicircular elements, see
\cite{VDN,NS}.

All polynomials, coefficient dimensions and entry laws are fixed as
$N\to\infty$. Superscript $\mathsf T$ denotes ordinary transpose
and $*$ denotes complex adjoint. A $*$-polynomial may involve the
variables and their formal adjoints, with arbitrary complex coefficients. We write
$\Tr$ for the ordinary trace and $\tr_d=d^{-1}\Tr$ for the normalized
trace. For a linear map $A$, its least singular value means
\[
 \smin(A)=\inf_{\|v\|=1}\|Av\|,
\]
where the infimum is over its domain. Singular values are otherwise
ordered decreasingly and counted with multiplicity. Constants may
depend on the fixed polynomial, coefficients and laws. We impose no
coupling between the matrices at different sizes.

\subsection{Independent matrices and their adjoints}
We consider a centered entry law $\xi$ with unit absolute variance:
\begin{equation}\label{eq:law}
 \EE\xi=0,\qquad \EE|\xi|^2=1.
\end{equation}
The density assumption takes one of the following two forms:
\begin{enumerate}[label=(\Alph*)]
\item $\xi$ is real and has a density bounded by $M$ with respect to
Lebesgue measure on $\R$;
\item $\xi$ is complex and has a joint density bounded by $M$ with
respect to Lebesgue measure on $\C=\R^2$.
\end{enumerate}
In class (B), the real and imaginary parts may be dependent. In
particular, the quantity $\rho=\EE\xi^2$ need not vanish.
Let $\xi_{\ell,ij}$ be independent copies of $\xi$ and define
\[
 X_\ell=X_\ell^{(N)}
 =\frac1{\sqrt N}(\xi_{\ell,ij})_{i,j=1}^N,
 \qquad 1\le\ell\le n.
\]
\begin{theorem}[Brown measure convergence]\label{thm:brown}
Assume either class (A) or class (B), and suppose that
$\EE|\xi|^{2+\varepsilon}<\infty$ for some $\varepsilon>0$.
Let $P$ be a fixed scalar-coefficient noncommutative $*$-polynomial.
Then
\[
 \mu_{P(X,X^*)}\xrightarrow{\ \PP\ }\mu_{P(c,c^*)}
\]
weakly, where $c_1,\ldots,c_n$ is a standard free circular family.
For every fixed $z\ne P(0)$,
\[
 \frac1N\log|\det(P(X,X^*)-zI_N)|
 \xrightarrow{\ \PP\ }\tau(\log|P(c,c^*)-z|).
\]
Moreover, the logarithmic potentials converge locally in planar
$L^1$ in probability: for every compact $K\subset\C$,
\[
 \int_K\bigl|U_{P(X,X^*)}(z)-U_{P(c,c^*)}(z)\bigr|\,\dd^2z
 \xrightarrow{\ \PP\ }0.
\]
\end{theorem}

Thus the same conclusion holds for real entry laws and arbitrary
complex polynomial coefficients, and for planar laws with arbitrary
pseudovariance $\rho$. The least singular value theorem used to prove
this result has a weaker moment assumption. To state it, fix
$A_\ell,B_\ell\in M_m(\C)$ and consider
\begin{equation}\label{eq:original-pencil}
 T_N=I_{mN}+\sum_{\ell=1}^n
       (A_\ell\otimes X_\ell+B_\ell\otimes X_\ell^*).
\end{equation}
We shall work on the Hilbert--Schmidt event
\begin{equation}\label{eq:Omega}
 \Omega_N=\left\{\max_\ell N^{-1}\norm{X_\ell}_{\HS}^2\le2\right\},
 \qquad q_N=\PP(\Omega_N^c),\qquad
 \omega=\begin{cases}1,&\text{class (A)},\\
                     1/2,&\text{class (B)}.\end{cases}
\end{equation}
Since the entries have a second moment, the weak law of large numbers
gives $q_N\to0$.

\begin{theorem}[Least singular values of tensor pencils]\label{thm:pencil}
Under \eqref{eq:law} and either density assumption (A) or (B),
there are $C,K,c>0$ and $N_0\in\N$, depending on the fixed
coefficients and law, such that
\begin{equation}\label{eq:pencil}
 \PP\{\smin(T_N)\le t\}
 \le CN^K t^{\omega/r_{\mathrm{comp}}}+Ce^{-cN}+q_N,
 \qquad N\ge N_0,\quad t>0.
\end{equation}
Here $r_{\mathrm{comp}}\le m$ is the length of any complex composition series of
$\C^m$ for the algebra generated by $I,A_\ell,B_\ell$.
For a nonzero irreducible coefficient family one may take $r_{\mathrm{comp}}=1$;
one may always replace $r_{\mathrm{comp}}$ by $m$.
If $\EE|\xi|^{2+\varepsilon}<\infty$, then
$q_N\le CN^{-\min(\varepsilon,2)}$.
For a real Gaussian law, or a centered complex Gaussian law with
positive definite real covariance, one may take $\omega=1$ and
replace the two additive errors by $Ce^{-cN}$.
\end{theorem}

The estimate is uniform in $t$, with constants allowed to depend on
the fixed coefficient family. The denominator $r_{\mathrm{comp}}$ arises from
reconstructing an upper triangular pencil from its irreducible
diagonal blocks.

Linearizing $P-zI$ gives a tensor pencil to which
Theorem~\ref{thm:pencil} applies. Its constant matrix is invertible
for $z\ne P(0)$.

\begin{corollary}[Least singular values of polynomial shifts]\label{cor:poly}
Let $P$ be a fixed scalar-coefficient $*$-polynomial, and assume
\eqref{eq:law} and either (A) or (B). For each fixed $z\ne P(0)$,
there are $C_z,a_z,b_z,c_z>0$ and $N_0(z)$ such that
\[
 \PP\{\smin(P(X,X^*)-zI_N)\le t\}
 \le C_zN^{a_z}t^{b_z}+C_ze^{-c_zN}+q_N,
 \qquad N\ge N_0(z),\quad t>0.
\]
In particular, for some fixed $B_z<\infty$,
\[
 \PP\{\smin(P(X,X^*)-zI_N)\ge N^{-B_z}\}\longrightarrow1.
\]
\end{corollary}

The constants may depend on $z$. The lower cutoff for each fixed
$z\ne P(0)$ suffices for Hermitization: the logarithmic potentials
are integrated over $z$, and the proof does not require uniformity in $z$.

\subsection{Wigner matrices}
For each $1\le\ell\le n$, let $H_\ell=H_\ell^{(N)}$ be Hermitian
with entries
\[
 (H_\ell)_{ij}=N^{-1/2}\zeta_{\ell,ij}\quad(i<j),\qquad
 (H_\ell)_{ji}=N^{-1/2}\overline{\zeta_{\ell,ij}},\qquad
 (H_\ell)_{ii}=N^{-1/2}d_{\ell,i}.
\]
All upper-triangle and diagonal atoms are independent, also across
$\ell$. For each fixed $\ell$, the off-diagonal atoms are copies of
a fixed centered law with absolute variance one. This law is either
real with a bounded one-dimensional density or complex with a bounded
joint planar density. Different input matrices may have different laws,
and both cases may occur in the same family. The real and imaginary parts
of a planar atom may be dependent. The diagonal atoms of each matrix are copies
of a fixed centered real law with bounded density and finite second
moment.

\begin{theorem}[Brown measure convergence for Wigner polynomials]\label{thm:wigner-brown}
Assume the preceding Wigner model, and suppose that the off-diagonal
laws have a finite $(2+\varepsilon)$ moment for some
$\varepsilon>0$. Let $P$ be a fixed scalar-coefficient polynomial
in the Hermitian variables. Then
\[
 \mu_{P(H_1,\ldots,H_n)}\xrightarrow{\ \PP\ }
 \mu_{P(s_1,\ldots,s_n)}
\]
weakly, where $s_1,\ldots,s_n$ is a standard free semicircular family.
For each fixed $z\ne P(0)$,
\[
 \frac1N\log|\det(P(H)-zI_N)|
 \xrightarrow{\ \PP\ }\tau(\log|P(s)-z|).
\]
The logarithmic potentials converge locally in planar $L^1$ in
probability. Only a second moment is required of the diagonal laws.
\end{theorem}

For Wigner matrices, we prove the least singular value bound using
their independent upper-triangle atoms.
For fixed $C_\ell\in M_m(\C)$, set
\[
 T_N^{\mathrm W}=I_{mN}+\sum_{\ell=1}^n C_\ell\otimes H_\ell,
\]
and define
\begin{equation}\label{eq:wigner-Omega}
 \Omega_N^{\mathrm W}=\left\{\max_{\ell\le n}
 N^{-1}\|H_\ell\|_{\HS}^2\le2\right\},
 \qquad q_N^{\mathrm W}=\PP((\Omega_N^{\mathrm W})^c).
\end{equation}
Let $\omega=1$ when all off-diagonal entry laws are real, and
$\omega=1/2$ when at least one has a planar density.

\begin{theorem}[Least singular values of Wigner pencils]\label{thm:wigner-pencil}
Under the stated density and second-moment assumptions there are
$C,K,c>0$ and $N_0\in\N$ such that
\[
 \PP\{\smin(T_N^{\mathrm W})\le t\}
 \le CN^K t^{\omega/r_{\mathrm{comp}}}+Ce^{-cN}+q_N^{\mathrm W},
 \qquad N\ge N_0,\quad t>0.
\]
Here $r_{\mathrm{comp}}\le m$ is the length of any complex composition series of
$\C^m$ for the algebra generated by $I,C_1,\ldots,C_n$, and
$q_N^{\mathrm W}\to0$.
If all atoms are Gaussian, one may take $\omega=1$ and replace
the two additive errors by $Ce^{-cN}$.
\end{theorem}

As in Corollary~\ref{cor:poly}, this theorem gives a polynomial
least singular value cutoff for $P(H)-zI_N$ at every fixed
$z\ne P(0)$. The proof is given in Section~\ref{sec:wigner}.
The single excluded shift has no effect on either weak eigenvalue
convergence or local spatial $L^1$ convergence.

\subsection{Main ideas of the proof}\label{intro:proof-ideas}
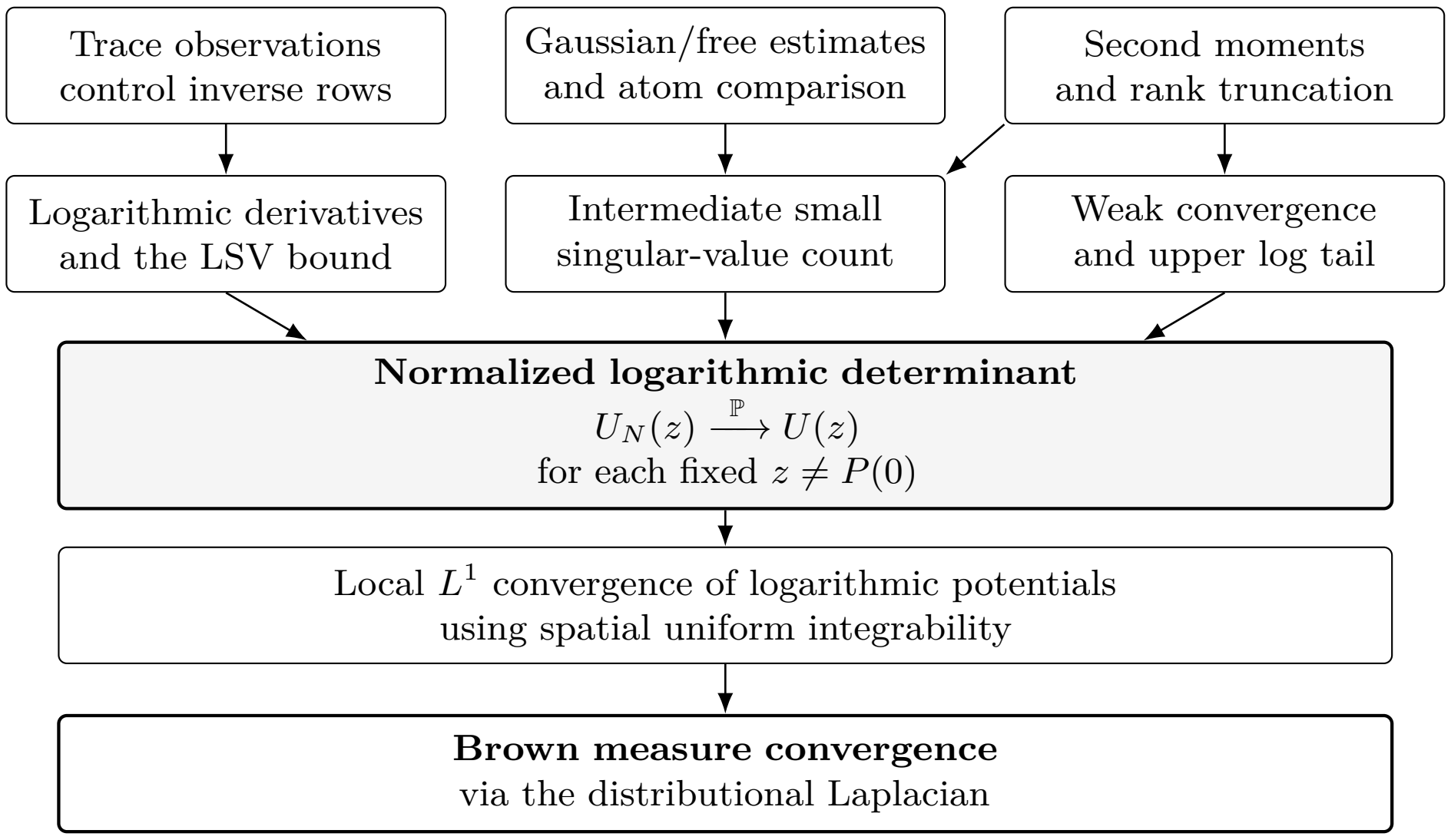
\begin{figure}[H]
\centering
\begin{tikzpicture}[
  box/.style={draw,rounded corners=2pt,align=center,
    font=\footnotesize,text width=3.5cm,minimum height=1.0cm,inner sep=4pt},
  wide/.style={box,text width=11.2cm},
  arrow/.style={-{Latex[length=1.8mm]},line width=.5pt}
]
\node[box] (algebra) at (-4.3,0)
  {Trace observations\\control inverse rows};
\node[box] (comparison) at (0,0)
  {Gaussian/free estimates\\and atom comparison};
\node[box] (moments) at (4.3,0)
  {Second moments\\and rank truncation};
\node[box] (lsv) at (-4.3,-1.45)
  {Logarithmic derivatives\\and the LSV bound};
\node[box] (count) at (0,-1.45)
  {Intermediate small\\singular-value count};
\node[box] (weak) at (4.3,-1.45)
  {Weak convergence\\and upper log tail};
\node[wide,line width=.8pt,fill=black!4] (logdet) at (0,-3.1)
  {\textbf{Normalized logarithmic determinant}\\
   \mbox{$U_N(z)\xrightarrow{\ \PP\ }U(z)$}\\
   \mbox{for each fixed $z\ne P(0)$}};
\node[wide] (local) at (0,-4.65)
  {Local $L^1$ convergence of logarithmic potentials\\
   using spatial uniform integrability};
\node[wide,line width=.8pt] (brown) at (0,-6.1)
  {\textbf{Brown measure convergence}\\
   via the distributional Laplacian};
\draw[arrow] (algebra) -- (lsv);
\draw[arrow] (comparison) -- (count);
\draw[arrow] (moments) -- (weak);
\draw[arrow] (moments.south west) -- (count.north east);
\draw[arrow] (lsv.south) -- ($(logdet.north)+(-3.6cm,0)$);
\draw[arrow] (count.south) -- (logdet.north);
\draw[arrow] (weak.south) -- ($(logdet.north)+(3.6cm,0)$);
\draw[arrow] (logdet) -- (local);
\draw[arrow] (local) -- (brown);
\end{tikzpicture}
\caption{From singular-value estimates to logarithmic potentials and Brown measures.}
\label{fig:proof-flow}
\end{figure}

The proof is organized around the normalized logarithmic determinant.
Write $P_N=P(X,X^*)$ or $P_N=P(H)$. Its logarithmic potential is
\[
 U_N(z)=\frac1N\log|\det(P_N-zI_N)|
       =\frac1N\sum_{j=1}^N\log s_j(P_N-zI_N),
 \qquad \mu_{P_N}=\frac1{2\pi}\Delta U_N.
\]
Thus the limiting eigenvalue measure can be identified through local
$L^1$ convergence of $U_N$ to the logarithmic potential $U$ of the free
polynomial. To justify the logarithmic integral, weak singular-value convergence
must be supplemented by estimates for very small and very large
singular values.
Figure~\ref{fig:proof-flow} shows how the estimates address these tails.

\begin{figure}[H]
\centering
\begin{tikzpicture}[
  calcbox/.style={draw,rounded corners=2pt,align=center,
    font=\footnotesize,inner sep=6pt},
  calcarrow/.style={-{Latex[length=1.8mm]},line width=.5pt}
]
\node[calcbox,text width=3.1cm,minimum height=2.2cm] (atom) at (-4.55,0)
  {};
\node[font=\footnotesize] at (-4.55,.78)
  {$\partial_uT_N$, $i\ne j$};
\begin{scope}[shift={(-4.55,-.42)},font=\footnotesize]
  \fill[black!9] (0,0) rectangle (.68,.58);
  \fill[black!9] (-.68,-.58) rectangle (0,0);
  \draw (-.68,-.58) rectangle (.68,.58);
  \draw (0,-.58) -- (0,.58);
  \draw (-.68,0) -- (.68,0);
  \node at (-.34,.29) {$0$};
  \node at (.34,.29) {$A$};
  \node at (-.34,-.29) {$B$};
  \node at (.34,-.29) {$0$};
  \node at (-.34,.76) {$i$};
  \node at (.34,.76) {$j$};
  \node at (-.87,.29) {$i$};
  \node at (-.87,-.29) {$j$};
  \node at (-1.24,0) {$\frac1{\sqrt N}$};
\end{scope}
\node[calcbox,text width=8.2cm,minimum height=2.2cm] (jacobi) at (1.9,0)
  {\textbf{Jacobi's formula} --- real coordinate $u=u_{\ell,ij,h}$\\[3pt]
   \mbox{$\displaystyle
     \frac{\partial_uF}{F}=\Tr(R\,\partial_uT_N)
     =\frac{\Tr(R_{ji}A)+\Tr(R_{ij}B)}{\sqrt N}$}\\[3pt]
   \mbox{$F=\det T_N,\quad R=T_N^{-1},\quad
            A=A_{\ell,h},\ B=B_{\ell,h}$}};
\draw[calcarrow] (atom.east) -- (jacobi.west);
\node[calcbox,text width=11.6cm,minimum height=1.1cm,fill=black!4,
      anchor=north] (inverse) at ($(jacobi.south)+(-1.9cm,-.55cm)$)
  {\textbf{Sum the inverse-row trace frame over all rows}\\[3pt]
   \mbox{$\displaystyle
     \|R\|_{\HS}^{2}
      \le N^{\beta}\sum_{\iota}
           \left|\frac{\partial_\iota F}{F}\right|^{2}
      =N^{\beta}\frac{\|\nabla F\|_2^{2}}{|F|^{2}}
      \qquad\text{on }\mathcal E_N$}};
\draw[calcarrow] (jacobi.south) -- (inverse.north);
\node[calcbox,text width=11.6cm,minimum height=1.1cm,anchor=north]
  (small) at ($(inverse.south)+(0,-.5cm)$)
  {\textbf{A small least singular value forces a relative determinant event}\\[3pt]
   \mbox{$\displaystyle
     \{\smin(T_N)\le t\}\cap\mathcal E_N
       \subseteq
     \{|F|\le N^{\beta/2}t\,\|\nabla F\|_2\}$}\\[3pt]
   Relative polynomial anti-concentration bounds the event on the right.};
\draw[calcarrow] (inverse) -- (small);
\end{tikzpicture}
\caption{Logarithmic derivatives convert inverse traces into determinant
anti-concentration. Only the affected spatial indices are shown.
Here $\iota$ ranges over the original real coordinates, and
$\mathcal E_N\subset\Omega_N$ is the frame event with exponent $\beta$
for a nonzero irreducible coefficient family. The zero set of $F$ is omitted.
For $i=j$, the two perturbed blocks merge into $A+B$.}
\label{fig:log-derivative}
\end{figure}
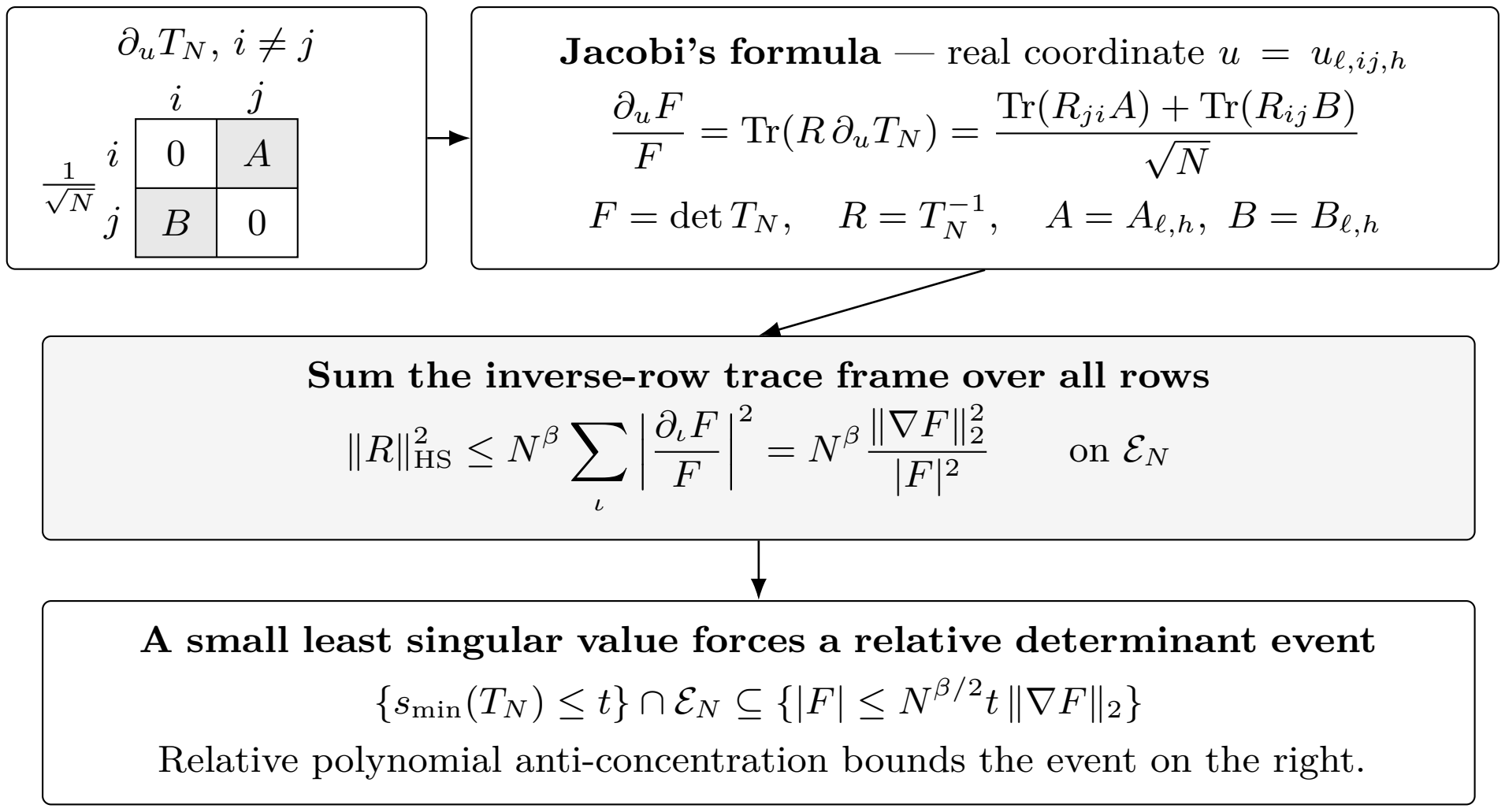

We follow the polynomial linearization and irreducible reduction of
\cite[Sections~2 and~3.1]{Han}. Its Gaussian and non-Gaussian results
treat polynomials without adjoints. On each nontrivial factor, we must show that traces
against the coefficients control the size of an inverse row.
Section~\ref{sec:algebra} constructs the required trace observations;
Section~\ref{sec:abstract} adapts the dimension-surplus covering
argument of \cite{Han} to the shared real atom coordinates. The inverse
equations then give a bound for every inverse row. Summing these bounds
and applying Jacobi's formula yields, for $F=\det T_N$,
\[
 \|T_N^{-1}\|_{\HS}\le N^C\frac{\|\nabla F\|_2}{|F|}
\]
on $\Omega_N$ outside an event of probability $Ce^{-cN}$.
Figure~\ref{fig:log-derivative} displays the calculation behind this bound.

This is the main departure from the least singular value argument in
\cite[Sections~3.4--3.7 and~6.1]{Han}. There one conditions on all but
one block column, estimates the determinant of its fixed-size
projection, and applies the block negative second moment identity.
Adjoint and Hermitian couplings make different block columns share
atoms. Here the trace observations give derivatives of the full
determinant in the original real coordinates. Relative polynomial
anti-concentration then gives the least singular value bound in
Section~\ref{sec:closing}. The triangular reconstruction from
\cite[Section~3.1]{Han} combines the factors in
Section~\ref{sec:reconstruction}; Section~\ref{sec:wigner} treats the
coupling at the two endpoints of a Wigner edge.

To prove convergence of $U_N$, we combine the least singular value
cutoff with a count at intermediate scales. Gaussian/free comparison
and rank truncation give this count, and together the bounds control
the lower logarithmic tail.
Second moments control the upper tail. Section~\ref{sec:spectral}
combines these estimates with weak singular-value convergence to obtain
pointwise convergence of $U_N$. A uniform expected spatial $L^2$ bound
gives uniform integrability and local $L^1$ convergence; the
Laplacian then identifies the Brown measure. Thus the pencil determinant
is used through its logarithmic derivatives to control the inverse,
and the polynomial determinant through its normalized logarithm to
identify the eigenvalue measure.

\subsection{Relation to previous work}
The free limits appearing here come from Voiculescu's asymptotic
freeness theorem \cite{Voiculescu}. For independent GUE matrices,
Haagerup and Thorbj{\o}rnsen \cite{HT} also proved convergence of
polynomial norms. Male \cite{Male} extended this framework to the
presence of independent matrices with a suitable strong limit, and
Anderson \cite{Anderson} proved polynomial norm convergence for
independent Wigner matrices under fourth-moment hypotheses.
Belinschi, Bercovici and Capitaine \cite{BBC} studied outliers of
self-adjoint polynomial deformations under spectral assumptions
on the input matrices. These results give norm convergence and
spectral exclusion in their respective settings. Convergence of
the eigenvalue measure for a nonnormal polynomial additionally
requires control of the singular logarithm.

Operator-valued subordination gives analytic descriptions of free
polynomial distributions \cite{BMS}. Linearizations and
realizations are developed in \cite{HMS}, while Hermitian reduction
provides a method for studying Brown measures of nonnormal free
polynomials \cite{BSS}. Cook, Guionnet
and Husson \cite{CGH} proved Brown measure convergence for
quadratic polynomials in independent complex Ginibre matrices.
The author's June 2026 paper \cite{Han} established this convergence
for polynomials of arbitrary degree without adjoints in both Ginibre
matrices and non-Gaussian iid matrices. The non-Gaussian convergence
holds under bounded planar density and finite moments of all orders;
see \cite[Theorem~1.4 and Section~6]{Han}. Its coefficient-module method
also gives least singular value estimates for linearized tensor pencils
in both settings, with the non-Gaussian estimate proved in
\cite[Theorem~6.3]{Han}.

Collins, Kitagawa, Miyagawa, Yuan and Zhong \cite{Collins} subsequently proved
weak Brown measure convergence for arbitrary GUE polynomials,
equivalently for arbitrary $*$-polynomials in independent complex
Ginibre matrices. Their least singular value theorem for Ginibre
pencils with structured means extends the Gaussian case of the centered
tensor-pencil estimates in \cite[Theorems~1.5 and~6.3]{Han}, following
its coefficient-space and restricted-adjoint argument; see \cite[Sections~1.3, 3.3 and~4]{Collins}.
They combine this extension with an auxiliary noisy ensemble,
a second-moment criterion for Brown measures, Gaussian rotation
and convexity. They also established convergence
of the spectral set and spectral radius, as well as a free support
theorem.

The polynomial convergence results in \cite[Theorems~1.1 and~1.4]{Collins}
cover Gaussian inputs only: independent GUE matrices, or equivalently independent complex
Ginibre matrices with adjoints. General Wigner matrices are explicitly
left untreated in \cite[Section~1.1]{Collins}. This restriction enters
the proof of \cite[Theorem~3.3]{Collins}: real rotations of two independent
GUE families preserve their joint law, making the relevant subharmonic
function independent of the real part of the rotation parameter and
convex in its imaginary part. For non-Gaussian entry laws, the same
rotation can change the joint law. The bounded densities and moment assumptions considered here do not
imply this invariance, so the rotation step does not yield the
non-Gaussian convergence theorem proved below.

The present proof develops the coefficient-module reduction and
restricted-adjoint argument of \cite{Han}. The doubled coefficient
algebra and its adjoint space give inverse-row trace observations that
control the full determinant gradient. This replaces conditional block
exposure and treats adjoint and Hermitian couplings. The least singular
value proof uses neither Gaussian rotation nor an auxiliary Gaussian
perturbation, and applies to real entry laws under the stated moment
assumptions.

Our probability estimates are also related to the scalar
invertibility theory. Littlewood--Offord and small-ball methods
were developed for independent-entry matrices in
\cite{RVLO,RVRect}, and decoupling methods for symmetric matrices
in \cite{VSym}. Under bounded densities, Farrell and Vershynin
\cite{FV} obtained polynomial invertibility bounds for symmetric
matrices without a moment assumption, uniformly over deterministic
symmetric additions. Here an atom can enter several coefficient
blocks, and the resulting tensor pencil can be nonnormal. The
coefficient algebra and the trace estimate for its inverse blocks
allow us to treat this dependence.

The Hermitian pencils are also related to the Kronecker
random matrices studied in \cite{AEKN}.
The spectral comparison uses the GUE resolvent estimate of Anderson, Guionnet and Zeitouni
\cite[Lemma~5.5.4]{AGZ} and the non-atomicity and algebraicity
results of Shlyakhtenko and Skoufranis
\cite[Theorem~1.1 and Corollary~1.2]{SS}. The former estimates the
Gaussian resolvents; the latter gives the small-singular-value bound
for the free limit.

\subsection{Real coordinates}
We treat classes (A) and (B) in a common notation. Set
$\mathsf d=1$ in class (A) and $\mathsf d=2$ in
class (B), and write
\begin{equation}\label{eq:block-atoms}
 u_{\ell,ij}=(u_{\ell,ij,h})_{h=1}^{\mathsf d},\qquad
 \xi_{\ell,ij}=\sum_{h=1}^{\mathsf d}\kappa_hu_{\ell,ij,h},
 \qquad\kappa_1=1,\quad\kappa_2=i.
\end{equation}
Here $u=\xi$ in class (A) and $u=(\Re\xi,\Im\xi)$ in
class (B). These independent blocks have bounded joint densities
and $\EE\|u\|^2=1$.

Define the real matrices
$Y_{\ell,h}=N^{-1/2}(u_{\ell,ij,h})_{ij}$ and the expanded
coefficients
\begin{equation}\label{eq:expanded-coeff}
 A_{\ell,h}=\kappa_hA_\ell,\qquad
 B_{\ell,h}=\overline{\kappa_h}B_\ell.
\end{equation}
Then \eqref{eq:original-pencil} becomes
\begin{equation}\label{eq:real-coordinate-pencil}
 T_N=I+\sum_{\ell,h}
 (A_{\ell,h}\otimes Y_{\ell,h}
       +B_{\ell,h}\otimes Y_{\ell,h}^{\mathsf T}).
\end{equation}
In particular, the imaginary coordinate contributes
$iA_\ell\otimes Y_{\ell,2}-iB_\ell\otimes Y_{\ell,2}^{\mathsf T}$.
The expanded coefficients generate the same complex algebra as
$A_\ell,B_\ell$: the $h=1$ coefficients already contain this
family, and the others are scalar multiples.

The real covariance of a planar atom is
\begin{equation}\label{eq:Sigma}
 \Sigma_\rho=\frac12
 \begin{pmatrix}1+\Re\rho&\Im\rho\\
                 \Im\rho&1-\Re\rho\end{pmatrix}.
\end{equation}
The bound $|\rho|\le1$ follows from \eqref{eq:law}. In class (B),
$\Sigma_\rho$ is positive definite and hence $|\rho|<1$.
Indeed, a zero-variance real linear form would force the law to be supported
on a line, contradicting its planar density. Class (A) corresponds to
$\rho=1$ and a one-dimensional real atom.

In class (A), $X_\ell^*=X_\ell^{\mathsf T}$.
For a complex matrix space viewed as a real Hilbert space we use
\[
 \langle U,V\rangle_\R=\Re\Tr(U^*V).
\]
Spatial blocks are taken after the fixed unitary tensor flip that places
the spatial indices first. In particular, the spatial blocks
of an $mN\times mN$ tensor pencil are $m\times m$ matrices.

\subsection{Organization}
Section~\ref{sec:linearization} identifies the trace observations for
linearized pencils. Section~\ref{sec:abstract} proves the probability
estimate that allows them to control inverse rows.
Section~\ref{sec:closing} passes from these rows to the determinant
gradient and the least singular value. Sections~\ref{sec:reconstruction}
and~\ref{sec:wigner} extend the bound to arbitrary coefficients and
Wigner inputs. Section~\ref{sec:spectral} controls the logarithmic tails
and proves Brown measure convergence.

\smallskip
\noindent\textit{AI assistance.}
The author used GPT-6 Astra and GPT-6.1 Sol to assist with algebraic
questions and some technical aspects of the $(2+\varepsilon)$-moment
arguments.

\section{Linearization and trace observations}\label{sec:linearization}\label{sec:algebra}

We first linearize the polynomial, following \cite[Section~2]{Han}.
The new task is to identify traces that control inverse rows when an
atom enters two transposed positions. A doubled coefficient algebra
expresses these paired traces. The resulting trace kernel satisfies
the invariant-subspace condition required in Section~\ref{sec:abstract}.

\subsection{Linearization of polynomial shifts}

We follow the reduction from polynomial shifts to normalized tensor
pencils in \cite[Section~2]{Han}. The word-chain construction below
also gives polynomial elimination matrices; see \cite{Anderson,BMS,HMS}
for the general theory. A matrix-valued $*$-polynomial has a fixed number of
blocks, each of which is a fixed scalar-coefficient noncommutative
$*$-polynomial in $Z,Z^*$. Evaluation at $N\times N$ matrices therefore
gives $N\times N$ blocks. The block number, degrees and coefficients do
not depend on $N$.

\begin{lemma}[Explicit word-chain linearization]\label{lem:linearization}
For every nonconstant $P$ with $c_0=P(0)$ there is a fixed-size
affine $*$-pencil
\[
 \mathcal L_z(Z)=K_z\otimes I+
   \sum_\ell(U_\ell\otimes Z_\ell+V_\ell\otimes Z_\ell^*)
\]
such that $K_z$ is invertible for $z\ne c_0$; the upper-left inverse
block of $\mathcal L_z$ is $(P(Z,Z^*)-zI)^{-1}$ whenever it exists;
and there are invertible matrix-valued $*$-polynomials $\mathcal V_z,\mathcal W_z$ with
\[
 \mathcal V_z\mathcal L_z\mathcal W_z=\diag(P(Z,Z^*)-zI,I).
\]
Both elimination matrices and their inverses are fixed matrix-valued
$*$-polynomials, uniformly bounded when input norms are bounded.
These are algebraic identities for every tuple of complex matrices
$Z$ and thus apply to both input classes with $Z=X$.
\end{lemma}
\begin{proof}
Write $P=c_0+\sum_w a_ww$, where the sum contains only
nonempty words with nonzero coefficients. For
$w=Z_{i_1}^{\varepsilon_1}\cdots Z_{i_{d_w}}^{\varepsilon_{d_w}}$,
with each $\varepsilon_j$ either $1$ or $*$, let $S_w$ have size
$d_w+1$ and first superdiagonal
\[
 (S_w)_{j,j+1}=Z_{i_j}^{\varepsilon_j},\qquad 1\le j\le d_w,
\]
with all other blocks zero. Nilpotence gives
\[
 (I-S_w)^{-1}=\sum_{k=0}^{d_w}S_w^k,
 \qquad ((I-S_w)^{-1})_{1,d_w+1}=w.
\]
Set $Q=\bigoplus_w(I-S_w)$. Let $u$ be the block row with
$-a_wI$ in the first coordinate of each word chain and zero
elsewhere, and let $v$ be the block column with $I$ in the last
coordinate of each chain. Then
\[
 uQ^{-1}v=-\sum_w a_ww=c_0-P.
\]
In the following matrices scalar constants are tensored with the
identity of the input dimension. Define
\[
 \mathcal L_z=
 \begin{pmatrix}c_0-z&u\\v&Q\end{pmatrix},\quad
 \mathcal V_z=
 \begin{pmatrix}I&0\\0&Q^{-1}\end{pmatrix}
 \begin{pmatrix}I&-uQ^{-1}\\0&I\end{pmatrix},\quad
 \mathcal W_z=\begin{pmatrix}I&0\\-Q^{-1}v&I\end{pmatrix}.
\]
Block multiplication, using $uQ^{-1}v=c_0-P$, gives
\[
 \mathcal L_z\mathcal W_z=
 \begin{pmatrix}P-zI&u\\0&Q\end{pmatrix},\qquad
 \mathcal V_z\mathcal L_z\mathcal W_z=\begin{pmatrix}P-zI&0\\0&I\end{pmatrix}.
\]
Whenever $P-zI$ is invertible, this identity gives
\[
 \mathcal L_z^{-1}
 =\mathcal W_z\begin{pmatrix}(P-zI)^{-1}&0\\0&I\end{pmatrix}\mathcal V_z.
\]
The first block row of $\mathcal W_z$ is $(I,0)$ and the first block
column of $\mathcal V_z$ is $(I,0)^{\mathsf T}$, which gives the
upper-left block of the inverse.

The finite geometric series for $Q^{-1}$ shows that $\mathcal V_z$ and
$\mathcal W_z$ are matrix-valued $*$-polynomials. Their unitriangular factors
invert by changing the off-diagonal sign, while
$\diag(I,Q^{-1})^{-1}=\diag(I,Q)$. Their inverses are also
matrix-valued $*$-polynomials, so all four matrices are uniformly bounded
when the input norms are bounded. At zero substitution
$Q=I$ and $uv=0$, since the first and last coordinates of each
nonempty word chain are distinct. The scalar Schur complement
then gives $\det K_z=c_0-z$, so $K_z$ is invertible exactly
when $z\ne c_0$.
\end{proof}

\begin{proof}[Proof of Corollary~\ref{cor:poly}]
For constant $P=c_0$, the least singular value equals the fixed
positive number $|c_0-z|$, so the assertion is immediate after
choosing the constant in the bound. Assume henceforth that $P$
is nonconstant.
Normalize by $T_{z,N}=(K_z^{-1}\otimes I)\mathcal L_z(X)$.
This is a pencil of the form in Theorem~\ref{thm:pencil}, with
fixed complex coefficients depending on $z$. Since
\[
 \norm{(P-zI)^{-1}}\le\norm{\mathcal L_z^{-1}}
 \le\norm{T_{z,N}^{-1}}\norm{K_z^{-1}},
\]
the desired bound follows.
\end{proof}

\begin{remark}[A singular shift in the real skew-symmetric model]
The exceptional shift cannot be removed uniformly over all
$*$-polynomials in the real model. For $P(x,x^*)=x-x^*$,
one has $P(0)=0$, and $P(X,X^{\mathsf T})$ is real
skew-symmetric. In every odd dimension its determinant is zero.
Excluding this scalar shift does not affect the almost-everywhere
argument in the plane for Brown measure convergence.
\end{remark}

\subsection{The doubled coefficient algebra}

Assume that $(A_\ell,B_\ell)$ is nonzero and complex irreducible.
The determinant derivatives pair traces at spatial indices $(i,j)$
and $(j,i)$. We adjoin a transpose copy of the pencil to express
these pairs as traces in one coefficient algebra. Its adjoint
space will be used in Theorem~\ref{thm:dual}.

By \eqref{eq:expanded-coeff}, the expanded tuple generates the same
algebra as the original tuple, so it is still nonzero and irreducible.
For every $(\ell,h)$ define
\begin{equation}\label{eq:double-coeff}
 C_{\ell,h,1}=
 \begin{pmatrix}A_{\ell,h}&0\\0&B_{\ell,h}^{\mathsf T}\end{pmatrix},
 \qquad
 C_{\ell,h,2}=
 \begin{pmatrix}B_{\ell,h}&0\\0&A_{\ell,h}^{\mathsf T}\end{pmatrix},
 \qquad
 \calB=\alg_{\C}(I,C_{\ell,h,1},C_{\ell,h,2}).
\end{equation}
The associated pencil, after the fixed coordinate ordering, is
\begin{equation}\label{eq:double-pencil}
 \widehat T_N=I+\sum_{\ell,h}
  (C_{\ell,h,1}\otimes Y_{\ell,h}+C_{\ell,h,2}\otimes Y_{\ell,h}^{\mathsf T})
 \simeq T_N\oplus T_N^{\mathsf T}.
\end{equation}
Indeed, its lower diagonal block is
\[
 I+\sum_{\ell,h}\bigl(B_{\ell,h}^{\mathsf T}\otimes Y_{\ell,h}
       +A_{\ell,h}^{\mathsf T}\otimes Y_{\ell,h}^{\mathsf T}\bigr)
 =T_N^{\mathsf T}.
\]
The transpose in \eqref{eq:double-pencil} is therefore the ordinary
transpose of $T_N$, with $X^*$ expressed in the expanded real
coordinates. The algebra $\calB$ is a complex
subalgebra of $M_m(\C)\oplus M_m(\C)$ and has complex dimension
at most $2m^2$.

\begin{lemma}[Trace pairing and common kernel]\label{lem:algebra}
The bilinear trace form $(b,c)\mapsto\Tr(bc)$ on $\calB$ is
nondegenerate. Moreover
$\bigcap_{\ell,h,t}\ker C_{\ell,h,t}=\{0\}$.
\end{lemma}
\begin{proof}
The first coordinate representation of $\calB$ on $\C^m$ is
irreducible by hypothesis. The second is also irreducible:
if a subspace is invariant under all $B_{\ell,h}^{\mathsf T},A_{\ell,h}^{\mathsf T}$,
its annihilator for the bilinear pairing $u^{\mathsf T}v$ is
invariant under all $B_{\ell,h},A_{\ell,h}$.
Thus the natural faithful $\calB$-module is a direct sum of two
irreducible modules. The radical annihilates each simple module,
so faithfulness forces the radical to vanish. Hence $\calB$ is
semisimple.

In a faithful representation of a finite-dimensional complex
semisimple algebra
$\bigoplus_jM_{d_j}(\C)$, the ambient trace is
$\sum_j r_j\Tr_{d_j}$ with positive integer multiplicities $r_j$.
If an element has a nonzero component in $M_{d_j}(\C)$, pairing
that component with a suitable matrix unit gives a nonzero trace.
Since $r_j$ is positive, the trace form on $\calB$ is
nondegenerate.

The common kernel of a nonzero irreducible family is zero: it is
an invariant subspace and cannot be the entire space. Apply this
to both coordinate families to obtain the last assertion.
\end{proof}

\subsection{Inverse blocks and trace observations}

The inverse blocks belong to the doubled algebra. Its traces pair the
two spatial positions changed by a single atom, as required for
determinant differentiation.

\begin{lemma}[Inverse blocks and the paired observations]\label{lem:inverse-blocks}
If $T_N$ is invertible and $R=T_N^{-1}$, every block
$\widehat R_{ij}$ of $\widehat R=\widehat T_N^{-1}$ belongs to
$\calB$, and
\begin{align}
 \widehat R_{ij}
 &=\diag(R_{ij},R_{ji}^{\mathsf T}),\label{eq:inverse-block}\\
 \Tr(\widehat R_{ji}C_{\ell,h,1})
 &=\Tr(R_{ji}A_{\ell,h})+\Tr(R_{ij}B_{\ell,h}),\label{eq:paired-A}\\
 \Tr(\widehat R_{ji}C_{\ell,h,2})
 &=\Tr(R_{ji}B_{\ell,h})+\Tr(R_{ij}A_{\ell,h}).\label{eq:paired-B}
\end{align}
\end{lemma}
\begin{proof}
The block formula follows from \eqref{eq:double-pencil}.
Cayley--Hamilton writes $\widehat T_N^{-1}$ as a polynomial in
$\widehat T_N$; its blocks therefore belong to the unital
algebra $\calB$. Multiplying the diagonal blocks and using
cyclicity and invariance of trace under transpose gives both trace
identities.
\end{proof}

\subsection{The trace kernel and invariant subspaces}

To control an inverse row by these traces, we must exclude profiles
concentrated in their common kernel. Let
$V=\calB^*:=\{b^*:b\in\calB\}$, with the real Hilbert--Schmidt
inner product, and define
\begin{equation}\label{eq:observation}
 \mathcal O(\Phi)=
    \bigl(\Tr(\Phi^*C_{\ell,h,t})\bigr)_{\ell,h,t},
 \qquad
 \calD=\ker\mathcal O\subset V.
\end{equation}
This map is real linear and complex conjugate linear, so its
kernel is a complex subspace. Adjoints of maps on $V$ below are
real Hilbert adjoints, including orthogonal compression to $V$
when the map is defined on the ambient matrix space; no closure
of $\calB$ under matrix adjoint is assumed. Also, for $\Phi=b^*\in V$,
\[
 C_{\ell,h,t}^*\Phi=(bC_{\ell,h,t})^*\in V.
\]
Thus left multiplication acts on $V$. The following lemma verifies
the structural hypothesis of Theorem~\ref{thm:dual}.

\begin{lemma}[Trace observability on the adjoint algebra]\label{lem:observable}
No nonzero real $E\subset\calD$ is invariant under all left
multiplications by $C_{\ell,h,t}^*$.
\end{lemma}
\begin{proof}
Suppose such a space $E$ exists. Since $\calD$ is complex,
$E_{\C}=E+iE$ is a nonzero complex subspace of $\calD$.
Invariance under the generators implies invariance under the
unital algebra $\calB^*$ acting on the left. Taking adjoints
reverses multiplication, so
\[
 J=\{\Phi^*:\Phi\in E_{\C}\}\subset\calB,
 \qquad Wb\in J\quad(W\in J,\ b\in\calB).
\]
Thus $J$ is a nonzero right ideal. Since every $Wb$ is the
adjoint of an element of $\calD$, for each generator we have
\[
 \Tr(C_{\ell,h,t}Wb)=\Tr(WbC_{\ell,h,t})=0
 \qquad\hbox{for every }b\in\calB.
\]
The element $C_{\ell,h,t}W$ belongs to $\calB$. Nondegeneracy
of the trace pairing therefore gives $C_{\ell,h,t}W=0$.
Every column of $W$ lies in the common kernel from
Lemma~\ref{lem:algebra}, so $W=0$. This contradicts $J\ne0$.
\end{proof}

\section{A probability estimate for the trace kernel}\label{sec:abstract}

The dimension surplus and multiscale covering argument come from
\cite[Theorem~3.3, Section~4 and Proposition~6.2]{Han}.
Here we adapt the probability estimates to whole real atom blocks:
coordinates within a block may be dependent, and the same atom
occurs in both row and column terms of the inverse equations.

We prove a polynomial least singular value bound for those equations
on the trace kernel. Taking the adjoint of an inverse row gives a vector
of coefficient blocks, called a profile, which solves every equation
except one. The preceding algebraic
argument supplies the structural condition that prevents such a
profile from concentrating in the trace kernel. The next section
uses this bound to control the full inverse through its traces.

For a fixed profile, anchor exposure leaves independent whole
blocks in different outputs. The dimension surplus allows summation
over the profile net. Spectral scales handle small singular values of
the coefficient map; projection handles small eigenvalues of the
profile Gram operator.

Let $V$ be a real Hilbert space of fixed dimension $p\ge1$, let
$\calD\subset V$ be a fixed subspace, and let
$L_{\ell,h},M_{\ell,h}:V\to V$ be fixed real linear maps.
Write $P_E$ for the orthogonal projection onto $E\subset V$ and
$\Gr(d,\calD)$ for the Grassmannian of $d$-dimensional real subspaces of $\calD$.
We use the independent blocks $u_{\ell,ij}\in\R^{\mathsf d}$
of \eqref{eq:block-atoms}, whose joint densities are bounded and
whose second moments are finite. The coordinates of a single block
may be dependent. Our structural assumption is
\begin{equation}\label{eq:no-invariant}
 \text{no nonzero real }E\subset\calD
 \text{ is invariant under every }L_{\ell,h},M_{\ell,h}.
\end{equation}
For a deleted output $a\in[N]$, define
\begin{equation}\label{eq:dual}
 (\frakD_a v)_k
 =v_k+\frac1{\sqrt N}\sum_{i,\ell,h}
  \bigl(u_{\ell,ik,h}L_{\ell,h}v_i+
        u_{\ell,ki,h}M_{\ell,h}v_i\bigr),\qquad k\ne a,
 \quad v\in\calD^N.
\end{equation}
The codomain is $V^{[N]\setminus\{a\}}$.
We equip profiles with $\norm v_2^2=\sum_i\norm{v_i}^2$,
and write $\widetilde{\frakD}_a$ for the same formula on $V^N$.
The Hilbert--Schmidt localization from the introduction gives
$\sum_h\norm{Y_{\ell,h}}_{\HS}^2=\norm{X_\ell}_{\HS}^2\le2N$
on $\Omega_N$. Summing the fixed coefficient norms and using
$\norm{Y_{\ell,h}},\norm{Y_{\ell,h}^{\mathsf T}}\le\sqrt{2N}$,
we obtain, simultaneously for all $a$,
\begin{equation}\label{eq:dual-norm}
 \norm{\widetilde{\frakD}_a}\le C_0\sqrt N.
\end{equation}

\begin{theorem}[Block-density dual estimate]\label{thm:dual}
If $\calD\ne0$, there are $E_*<\infty$ and $N_0$ such that, for $N\ge N_0$,
\begin{equation}\label{eq:dual-tail}
 \PP\left(\Omega_N\cap
 \{\exists a\in[N]:\smin(\frakD_a)\le N^{-E_*}\}\right)
 \le Ce^{-cN}.
\end{equation}
The constants depend only on the fixed spaces, maps, block dimension,
and block law. All spaces and dimensions in the proof are real.
The second moment ensures $q_N\to0$; the conditional small-ball
estimate itself uses bounded densities and a fixed cutoff.
\end{theorem}

\subsection{The dimension surplus and the covering geometry}

To obtain the dimension surplus from \eqref{eq:no-invariant},
list the coefficient maps as
$\mathsf H_1,\ldots,\mathsf H_r$, where $r=2n\mathsf d$, and for a nonzero real
$E\subset\calD$ set
\[
 \Psi_E:E^r\longrightarrow V,\qquad
 \Psi_E(v_1,\ldots,v_r)=\sum_j\mathsf H_jv_j.
\]
The image of $\Psi_E$ contains all possible random terms produced
by profiles in $E^N$. Some of its singular directions may, however,
be too small for a useful small-ball estimate. For $\sigma>0$, let
$Y_\sigma(E)$ be the left singular space corresponding to singular
values at least $\sigma$. For $0<\tau<1$, let
$C_{\tau,\sigma}(E)$ be the right singular space of
$P_{Y_\sigma(E)^\perp}|_E$ corresponding to singular values at most
$\tau$. Thus $C_{\tau,\sigma}(E)$ consists of the directions in
$E$ that are close to this active image.

These two singular spaces have different roles. Projection onto the
active image is where bounded densities give a small-ball estimate.
Projection onto its orthogonal complement instead controls the
possible profile coordinates: if the random contribution there is
small, the profile coordinate must be close to the active image.
The second singular space describes the directions that remain
free in this restriction. We show that the active image has at
least one more dimension, and then construct the corresponding
covering net.

The following lemma is the real-space counterpart of
\cite[Proposition~4.1]{Han}. Its surplus of one dimension makes the
small-ball estimate dominate the profile net.

\begin{lemma}[Uniform surplus]\label{lem:surplus}
There are $\sigma_0,\tau\in(0,1)$ such that for all nonzero
$E\subset\calD$ and $0<\sigma\le\sigma_0$,
\[
 y:=\dim Y_\sigma(E)\ge
 c_E+1,\qquad c_E:=\dim C_{\tau,\sigma}(E).
\]
Also $\norm{\Psi_E}\ge\sigma_0$.
\end{lemma}
\begin{proof}
By \eqref{eq:no-invariant} and compactness of the finite union of
real Grassmannians,
\[
b:=\inf_{0\ne E\subset\calD}
       \sum_j\norm{P_{E^\perp}\mathsf H_jP_E}>0.
\]
Let $M_0=\sup_E\norm{\Psi_E}<\infty$, and write
$Y=Y_\sigma(E)$ and $C=C_{\tau,\sigma}(E)$.
For $x\in C$ we have
$\norm{P_Yx}^2\ge(1-\tau^2)\norm x^2$, so $P_Y:C\to Y$
is injective and $c_E\le y$. Suppose equality holds. Then every
$w\in Y$ can be written as $w=P_Yx$ for an $x\in C$, and
\[
 \dist(w,E)\le\norm{w-x}
 \le\tau\norm x
 \le\frac{\tau}{\sqrt{1-\tau^2}}\norm w.
\]
For a unit $v\in E$, decompose
$\mathsf H_jv=P_Y\mathsf H_jv+P_{Y^\perp}\mathsf H_jv$. The second term has norm at
most $\sigma$, because it belongs to the discarded part of
$\Psi_E$. Applying the preceding distance bound to the first
term gives
\[
 \dist(\mathsf H_jv,E)\le
 M_0\frac{\tau}{\sqrt{1-\tau^2}}+\sigma.
\]
Taking operator norms and summing over $j$, we obtain
\[
 b\le r\left(M_0\frac{\tau}{\sqrt{1-\tau^2}}+\sigma\right).
\]
For example, take
\[
 \tau=\min\left(\frac14,\frac{b}{8r(1+M_0)}\right),
 \qquad
 \sigma_0=\min\left(\frac14,\frac{b}{4r}\right).
\]
The displayed upper bound is then at most $b/2$, a contradiction.
Only a positive lower bound for $b$ is used.
The same argument, or $y\ge1$, gives $\norm{\Psi_E}\ge\sigma_0$.
\end{proof}

\begin{lemma}[Real tube net]\label{lem:tube}
Let $Y=Y_\sigma(E)$, and suppose $\dim E=d$ with $c_E$ as above.
For $r_{\mathrm{tube}}\ge0$, the set
$\{v\in E:\norm v\le2,\ \dist(v,Y)\le r_{\mathrm{tube}}\}$ has an $h$-net of size
\[
 (C/h)^{c_E}
 \left(\frac{C\max\{r_{\mathrm{tube}}/\tau,h\}}{h}\right)^{d-c_E},
 \qquad 0<h\le1.
\]
The centers may be chosen in the set.
\end{lemma}
\begin{proof}
Write $T=P_{Y^\perp}|_E$ and decompose $v=v_0+v_1$, with
$v_0\in C_{\tau,\sigma}(E)$ and
$v_1\in E\cap C_{\tau,\sigma}(E)^\perp$. The singular-vector
decomposition gives
\[
 \norm{Tv}^2=\norm{Tv_0}^2+\norm{Tv_1}^2,
 \qquad \norm{Tv_1}\ge\tau\norm{v_1}.
\]
Thus $\norm{v_0}\le2$ and $\norm{v_1}\le r_{\mathrm{tube}}/\tau$.
Taking Euclidean ball nets in these two orthogonal spaces proves
the bound. Finally, moving each net center whose ball meets the
set to a point of that intersection changes the mesh by only a
fixed factor.
\end{proof}

\subsection{A fixed profile: anchors and independent outputs}

The row and column terms share atoms, so independence must be
obtained before multiplying output probabilities. We expose all
blocks except those joining a small anchor set to the other outputs,
and discard the anchor outputs. Each unexposed atom then enters
only one remaining equation. The profile values at the anchors
must span the active directions. For $v\in E^N$,
define its real Gram operator by
\[
 \Gamma_v=\sum_i v_i\otimes_{\R}v_i.
\]
Here $x\otimes_{\R}x$ denotes the map
$y\mapsto\langle x,y\rangle_\R x$.
We choose $d=\dim E$ indices with the following lower bound on
their Gram operator.

\begin{lemma}[Anchor extraction]\label{lem:anchors}
If $1\le\dim E=d\le p$, $\gamma>0$, $\norm v_2^2\le C_1$, and
$\Gamma_v\succeq\gamma P_E$, there is an $I\subset[N]$ of size $d$
such that
\[
 \Gamma_{v,I}:=\sum_{i\in I}v_i\otimes_{\R}v_i
 \succeq c_{d,C_1}\gamma^dN^{-d}P_E.
\]
\end{lemma}
\begin{proof}
Write the profile as a real $d$ by $N$ matrix $B$ in an
orthonormal basis of $E$. By Cauchy--Binet,
\[
 \gamma^d\le\det(BB^{\mathsf T})
 =\sum_{|I|=d}\det(B_I)^2.
\]
Choose a maximal minor. Computing determinants on $E$, we obtain
$\det\Gamma_{v,I}\ge\gamma^d/\binom Nd$, whereas every
eigenvalue of $\Gamma_{v,I}$ is at most
$\Tr\Gamma_{v,I}\le C_1$. Its least eigenvalue is therefore at
least $C_1^{-(d-1)}\gamma^d/\binom Nd$, which proves the lemma.
\end{proof}

Bounded joint density alone does not bound the density of a
projection from an unbounded set. We therefore restrict the
unexposed blocks to a fixed ball before applying the volume
estimate below. A small proportion of outputs may fail this
restriction; the following event bounds that proportion
simultaneously for every possible anchor set.
Fix $\kappa=1/(8p)$ and choose a finite $R_0\ge1$ so large that
\[
\PP\{\norm u\le R_0\}^{2np}\ge1-\kappa/4 .
\]
Under the normalization $\EE\norm u^2=1$, the choice
$R_0=8p\sqrt n$ suffices by Chebyshev's inequality and a union
bound over $2np$ whole blocks. More generally, existence of a
fixed cutoff requires only that the block law be a probability
measure. There is an event $\mathcal A_N$ of probability at least $1-Ce^{-cN}$
such that, simultaneously for all $I\subset[N]$ with $1\le|I|\le p$,
at most $\kappa N$ indices $k\notin I$ fail
\begin{equation}\label{eq:good-row}
 \norm{u_{\ell,ik}}\le R_0,\quad\norm{u_{\ell,ki}}\le R_0
       \qquad (i\in I,\ \ell\in[n]).
\end{equation}
For a fixed $I$, the indicators in question involve disjoint sets
of atoms as $k$ varies. Each has failure probability at most
$\kappa/4$, so Chernoff's inequality bounds the probability of
more than $\kappa N$ failures by $Ce^{-cN}$. Taking a union over
at most $pN^p$ choices of $I$ gives the assertion.

\begin{lemma}[Bounded-support projection]\label{lem:projection}
If a random vector $Z\in\R^D$ has density at most $K$ and support
in a ball of fixed radius $R_{\mathrm{cut}}$, and $A:\R^D\to\R^y$ has all $y$ singular values
at least $a>0$, then, uniformly in $w$,
\[
 \PP\{\norm{AZ-w}\le\delta\}\le C_{D,K,R_{\mathrm{cut}}}(\delta/a)^y .
\]
\end{lemma}
\begin{proof}
In singular-vector coordinates, $\norm{AZ-w}\le\delta$
restricts the first $y$ coordinates to an ellipsoid of volume
at most $C(\delta/a)^y$. The other coordinates remain in a
fixed ball. The preimage therefore has volume at most
$C_{D,R_{\mathrm{cut}}}(\delta/a)^y$, and multiplication by the density bound
proves the claim.
\end{proof}

We can now combine the anchor lower bound with the cutoff. In the
conditioning below each atom remains a whole block, so dependence
between its real coordinates does not affect the argument.

\begin{lemma}[One deterministic profile]\label{lem:pointwise}
In Lemma~\ref{lem:anchors}, write $y=\dim Y_\sigma(E)$ and
$b_N(d)=\lfloor(1-\kappa)N\rfloor-d-1$.
For $0<\delta\le1$, $0<\sigma\le\sigma_0$, and sufficiently large $N$,
\begin{equation}\label{eq:pointwise}
 \PP\{\norm{\frakD_av}_2\le\delta,\ \mathcal A_N\}
 \le\min\left\{1,\
 \left(C\sigma^{-y}\gamma^{-dy/2}
       N^{(d+1)y/2}\delta^y\right)^{b_N(d)}\right\}.
\end{equation}
For fixed $C_1$, the constants and the lower threshold on $N$
are uniform over the deterministic profiles, their subspaces,
$a$, $\gamma$, and $\sigma$ satisfying these assumptions.
\end{lemma}
\begin{proof}
Choose $I$ from Lemma~\ref{lem:anchors} for the fixed profile,
and let $J=[N]\setminus(I\cup\{a\})$. First expose all atoms
except the whole directed blocks $u_{\ell,ik},u_{\ell,ki}$ with
$i\in I$ and $k\in J$. For each output $k\in J$, its
unexposed random part is
\[
 Z_k=\frac1{\sqrt N}\sum_{i\in I,\ell,h}
       (u_{\ell,ik,h}L_{\ell,h}v_i+
        u_{\ell,ki,h}M_{\ell,h}v_i).
\]
The other occurrence of an unexposed atom lies in an anchor output,
which is discarded. Consequently the $Z_k$ use disjoint atoms as
$k$ varies; all other terms are fixed translations.

Set $Y=Y_\sigma(E)$, and let $A_k$ be the real linear map from
the $2nd\mathsf d$ coordinates of the unexposed blocks to $P_YZ_k$.
Its Euclidean Gram operator is
\[
 \frac1N P_Y\sum_{\ell,h}
 (L_{\ell,h}\Gamma_{v,I}L_{\ell,h}^*+
  M_{\ell,h}\Gamma_{v,I}M_{\ell,h}^*)P_Y
 \succeq c\sigma^2\gamma^dN^{-d-1}I_Y .
\]
Indeed, the anchor lower bound replaces $\Gamma_{v,I}$ by
$c\gamma^dN^{-d}P_E$, and
$P_Y\Psi_E\Psi_E^*P_Y\succeq\sigma^2I_Y$ by definition of
$Y$. Hence every singular value of $A_k$ onto $Y$ is at least
$c\sigma\gamma^{d/2}N^{-(d+1)/2}$. This computation uses
Euclidean coordinates and real Hilbert adjoints; it requires no
assumption on the random covariance inside a block.

Condition on the membership pattern of the good rows
\eqref{eq:good-row} in $J$. The rows remain independent, and
within each good row the blocks are independent copies of
$u$ conditioned on $\norm u\le R_0$.
If $p_0=\PP(\norm u\le R_0)>0$, their product density is at most
$(M/p_0)^{2nd}$ on a ball of radius $\sqrt{2nd}\,R_0$ in
$\R^{2nd\mathsf d}$. The dimension is at most $2np\mathsf d$,
so both the density bound and the support radius are uniform in the profile and $N$.
Applying Lemma~\ref{lem:projection} to $A_k$ gives the one-output
factor
\[
 \varpi=C\sigma^{-y}\gamma^{-dy/2}N^{(d+1)y/2}\delta^y,
\]
uniformly in the exposed translations and in the good-row pattern.

On $\mathcal A_N$, this particular $I$ has at least $b_N(d)$
good outputs in $J$. We enlarge to that event, multiply the
one-output estimates, and then average over the membership pattern.
More explicitly, if $\mathcal F_I$ denotes the first exposure and
$g$ ranges over patterns with at least $b_N(d)$ good rows, the
conditional bound is
\[
 \sum_g\PP(g\mid\mathcal F_I)\,
       \min(1,\varpi)^{\#\mathrm{good}(g)}
 \le \min(1,\varpi)^{b_N(d)}.
\]
The pattern probabilities sum to at most one. Integrating over
$\mathcal F_I$ proves \eqref{eq:pointwise}, with
$\mathcal A_N$ used only through the preceding enlargement,
rather than as a conditioning event.
\end{proof}

The estimate for a fixed profile has now been integrated over the
exposed atoms. Consequently
each point of a deterministic net may use its own anchor set when
we sum these probabilities. We construct this net on $\Omega_N$,
where the operator norm controls approximation errors; the event
$\Omega_N$ does not enter any conditional law.

\subsection{Uniform estimates for regular profiles}

We first treat profiles whose Gram operator has a lower bound on
their span. Write
\[
 \mathcal R(E,\gamma)=
 \{v\in E^N:\norm v_2\le1,\ \Gamma_v\succeq\gamma P_E\}.
\]
Norms below one allow the later projection of a unit profile.
The Gram lower bound makes anchor extraction uniform. A spectral
gap for $\Psi_E$ bounds its inactive contribution, forcing the
profile coordinates into the tubes above. The following lemma
assumes both bounds; the next two subsections remove them.

\begin{lemma}[Real uniform net bound]\label{lem:regular}
There are fixed $C_*,c_*,K_0,L_0,N_*>0$ with the following property,
uniformly in all the scale parameters below.
For definiteness one may take
$K_0=L_0=2p^2+4p+2$.
Let $\mathcal F\subset\Gr(d,\calD)$ be the closed family for which
all singular values of $\Psi_E$ lie in $[0,\eta]\cup[\sigma,\infty)$,
where $0<\eta<\sigma\le\sigma_0$.
For $N\ge N_*$, if
\[
 C_*\sqrt N\,\eta\le\delta,\qquad
 0<\delta\le c_*\gamma\le c_*,
\]
then
\begin{align}
 &\PP\bigl(\Omega_N\cap\mathcal A_N\cap
 \{\exists E\in\mathcal F,\ v\in\mathcal R(E,\gamma):
                      \norm{\frakD_av}_2\le\delta\}\bigr)
 \notag\\
 &\qquad\le
 C_*^N N^{K_0N}\gamma^{-p^2N}\sigma^{-pN}
                   \delta^{N/2-L_0} .\label{eq:regular}
\end{align}
\end{lemma}
\begin{proof}
\emph{Approximating the subspace.}
Choose a deterministic net of $\mathcal F$ in projection norm at
mesh $\theta=c\delta/\sqrt N$, with centers in $\mathcal F$.
Orthogonal projections lie in a bounded subset of the real
$p$ by $p$ matrices, so at most $(C/\theta)^{p^2}$ centers are
needed. For a nearby center $E'$, define $v_i^0=P_{E'}v_i$.
Since the two subspaces have equal dimension and projection
distance at most $\theta$, we have
$P_{E'}P_EP_{E'}\succeq(1-\theta^2)P_{E'}$. It follows that
\[
 \norm{v^0}_2\le1,\quad \norm{v-v^0}_2\le\theta,\quad
 \Gamma_{v^0}\succeq\gamma(1-\theta^2)P_{E'}
                         \succeq(15\gamma/16)P_{E'}.
\]
On $\Omega_N$, \eqref{eq:dual-norm} gives
$\norm{\frakD_av^0}_2\le2\delta$.

\emph{Restricting and covering the profile coordinates.}
Fix $Y=Y_\sigma(E')$. For each $k\ne a$, write the random sum
in \eqref{eq:dual} as $\Psi_{E'}t_k$, where the components of
$t_k\in(E')^r$ are
$N^{-1/2}\sum_i u_{\ell,ik,h}v_i^0$ and
$N^{-1/2}\sum_i u_{\ell,ki,h}v_i^0$. Cauchy--Schwarz gives
\[
 \norm{t_k}^2
 \le\frac2N\sum_{\ell,i,j}\norm{u_{\ell,ij}}^2
                 \norm{v^0}_2^2
 \le4nN
 \qquad\text{on }\Omega_N.
\]
This energy bound uses all coordinates before projection and needs no
further exposure of atoms.
The spectral gap at the net center implies
$\norm{P_{Y^\perp}\Psi_{E'}}\le\eta$. Consequently,
\[
 \dist(v_k^0,Y)
 \le\norm{(\frakD_av^0)_k}+
                 \norm{P_{Y^\perp}\Psi_{E'}t_k}
 \le2\delta+C\sqrt N\eta\le C\delta.
\]
The small-image equation thus forces each of these profile
coordinates into a thin tube about the active image. This
restriction is stronger than a ball bound and reduces the number
of directions that must be covered at the small mesh. Since the
images are taken at the net centers, we do not need continuity of
$Y_\sigma(E)$.

Apply Lemma~\ref{lem:tube} at mesh $h=c\delta/N$ to these
$N-1$ coordinates, and use a ball net for the removed coordinate.
If $c'=\dim C_{\tau,\sigma}(E')$, covering the $c'$ tube directions requires
$(CN/\delta)^{c'}$ per output; the remaining directions have
radius $O(\delta)$ and require at most $(CN)^p$ net points. The product net
therefore has size at most
\begin{equation}\label{eq:net-size}
 (CN/\delta)^{c'(N-1)+d}(CN)^{p(N-1)} .
\end{equation}
Rounding incurs profile error at most $c\delta/\sqrt N$, so its
image error on $\Omega_N$ is at most $C_0c\delta$. The inequality
\[
 \norm{\Gamma_w-\Gamma_{v^0}}
 \le(\norm w_2+\norm{v^0}_2)\norm{w-v^0}_2
\]
shows, using $\delta\le c_*\gamma$ and taking the fixed mesh
constant small enough, that the rounded profile satisfies
$\norm w_2\le2$ and $\Gamma_w\succeq(\gamma/2)P_{E'}$.
Keep only net points with these two properties. Every original
profile in the event is represented by one of the selected points,
and that point satisfies $\norm{\frakD_aw}_2\le C\delta$.

\emph{Summing the pointwise estimates.}
The preceding inclusion holds on $\Omega_N$, so we may drop this event before
applying Lemma~\ref{lem:pointwise} to each deterministic net point.
Its factor $\delta^{yb_N(d)}$, combined with the profile net,
leaves the exponent
\[
 \begin{split}
 yb_N(d)-c'(N-1)-d
 &\ge ((1-\kappa)y-c')N-y(d+2)+c'-d\\
 &\ge (7/8)N-C_p ,
 \end{split}
\]
because $y\ge c'+1$, $y\le p$, and $\kappa=1/(8p)$.
The subspace net contributes the additional fixed power
$\delta^{-p^2}$, after which we may weaken
the exponent to $N/2-L_0$. The remaining factors are bounded by
$C_*^N N^{K_0N}\gamma^{-p^2N}\sigma^{-pN}$, using $d,y\le p$
and $0<\gamma,\sigma\le1$.
More explicitly, their total power of $N$ is at most
\[
 2pN+\frac{p(p+1)}2N+\frac{p^2}2
 \le \left(p^2+\frac{5p}2\right)N.
\]
The original power of $\delta$ is at least
$7N/8-2p^2-3p$. These two bounds allow the stated choices
$K_0=L_0=2p^2+4p+2$. We take $c_*$ smaller if necessary so that
the radius $C\delta$ used in Lemma~\ref{lem:pointwise} is at
most one. This proves \eqref{eq:regular}.
\end{proof}

\subsection{Removing the spectral gap assumption}

The preceding lemma cannot be applied at a single fixed cutoff
for $\Psi_E$: a singular value could lie just below that cutoff,
too small for the density estimate but too large to discard.
We instead choose finitely many decreasing scales. Since the
target dimension is fixed, one of their intervals is empty, giving
the required gap. The recursion below ensures that the inactive
contribution at the lower end is smaller than the small-image
radius chosen at the upper end.

\begin{lemma}\label{lem:all-regular}
For every fixed $D\ge1$ there is a finite integer $e(D)$ such that,
for each $1\le d\le\dim\calD$,
\[
 \PP\bigl(\Omega_N\cap\mathcal A_N\cap
 \{\exists E\in\Gr(d,\calD),\
 v\in\mathcal R(E,N^{-D}):
       \norm{\frakD_av}_2\le N^{-e(D)}\}\bigr)
 \le p e^{-2N\log N}
\]
for all sufficiently large $N$.
\end{lemma}
\begin{proof}
Set
\[
 A=4p+8,\qquad B=4(K_0+p^2D+4),\qquad
 f_D(t)=t^A N^{-B}.
\]
For $N\ge4L_0$, $\gamma=N^{-D}$ and $\delta=f_D(\sigma)$,
the bound \eqref{eq:regular} is at most
\[
 \bigl(C_*N^{K_0+p^2D}\sigma^{-p}\delta^{1/4}\bigr)^N
 =\bigl(C_*\sigma^2N^{-4}\bigr)^N
 \le e^{-2N\log N}.
\]
Set $t_0=\sigma_0$ and
$t_{j+1}=f_D(t_j)N^{-3}$ for $0\le j<p$.
Lemma~\ref{lem:surplus} gives at least one singular value of
$\Psi_E$ at or above $t_0$. At most $p-1$ other positive
singular values can therefore lie in the $p$ disjoint intervals
$[t_{j+1},t_j)$. One interval is empty. All singular values then
belong to $[0,t_{j+1}]\cup[t_j,\infty)$, so this subspace is in
the gap family of Lemma~\ref{lem:regular} with
$\sigma=t_j$, $\eta=t_{j+1}$, and $\delta_j=f_D(t_j)$.
The contribution of the inactive singular directions is bounded by
\[
 C_*\sqrt N\,t_{j+1}
 =C_*N^{-5/2}\delta_j\le\delta_j.
\]
Also $\delta_j\le c_*N^{-D}$ for large $N$, since $B>D$.

The recursion gives
\[
 t_j=\sigma_0^{A^j}
       N^{-(B+3)(A^j-1)/(A-1)}.
\]
Any integer $e(D)$ greater than
$A(B+3)(A^{p-1}-1)/(A-1)+B+2$ ensures
$N^{-e(D)}\le\min_{j<p}\delta_j$ for all sufficiently large $N$.
A union over the $p$ gap families proves the lemma.
\end{proof}

\subsection{Projection to profiles with a Gram lower bound}

The spectral gap has been removed, but the last lemma still assumes
a polynomial lower bound for the profile Gram operator on its
span. An arbitrary unit profile need not satisfy such a bound.
We diagonalize its Gram operator and discard the directions with
small eigenvalues, projecting every component onto the same
remaining subspace. The squared error is the sum of the discarded
eigenvalues. The thresholds below make this error small enough,
even after multiplication by the localized operator norm, to apply
the preceding estimate to the projected profile.

\begin{proof}[Proof of Theorem~\ref{thm:dual}]
Let $q=\dim\calD$. Choose integers
\[
 D_1=1,\qquad D_{d+1}\ge
 \max\{D_d+1,\,2e(D_d)+4\},\quad 1\le d<q.
\]
For a unit profile $v\in\calD^N$, diagonalize its real Gram operator.
Starting at dimension $q$, discard its smallest eigenvector whenever
the smallest remaining eigenvalue is below $N^{-D_d}$ at dimension
$d$. Project every profile component at each step, without
renormalizing. The process stops at some $d\ge1$, because the largest
original eigenvalue is at least $1/q>N^{-1}$.
Let $E$ be the remaining eigenspace and $v'$ the projected profile.
Then $v'\in\mathcal R(E,N^{-D_d})$. At a discarded dimension
$j>d$, the squared projection error is the discarded eigenvalue,
which is less than $N^{-D_j}$. These errors are orthogonal, so
\[
 \norm{v-v'}_2^2
 \le\sum_{j=d+1}^qN^{-D_j}
 \le qN^{-D_{d+1}}\qquad(d<q).
\]
There is no error when $d=q$. On $\Omega_N$, \eqref{eq:dual-norm}
and the choice of $D_{d+1}$ now give
\[
 \norm{\frakD_a(v-v')}_2
 \le C_0\sqrt{qN}\,N^{-D_{d+1}/2}
 \le\frac12N^{-e(D_d)}
\]
for all sufficiently large $N$. Choose
$E_*>\max_d e(D_d)+1$. If
$\norm{\frakD_av}_2\le N^{-E_*}$, the triangle inequality
then gives $\norm{\frakD_av'}_2\le N^{-e(D_d)}$.
Lemma~\ref{lem:all-regular} bounds this event, within
$\Omega_N\cap\mathcal A_N$, by $pe^{-2N\log N}$ for each
possible $d$ and $a$. Taking the union over them and adding
$\PP(\mathcal A_N^c)\le Ce^{-cN}$ once proves
\eqref{eq:dual-tail}.

The families with a spectral gap and the profile sets with a Gram
lower bound used above are compact and defined by continuous maps. Their
existential events are therefore measurable, as is the least
singular value, which is continuous in the entries.
\end{proof}

\section{From inverse rows to the least singular value}\label{sec:closing}

Assume that the coefficient tuple is nonzero and complex irreducible.
The trace observations of Section~\ref{sec:algebra} control an inverse
row because its equations prevent concentration in the kernel of
$\mathcal O$. We sum these row bounds to control the full inverse,
then identify the paired traces with logarithmic derivatives of the
full determinant. Relative anti-concentration gives the least singular
value bound.

\subsection{The inverse equations and a frame estimate}

On the invertibility set, write $\widehat R=\widehat T_N^{-1}$.
For $a\in[N]$, set $\Phi_i=\widehat R_{ai}^*\in V$.
The $(a,k)$ block of $\widehat R\widehat T_N=I$ vanishes when
$k\ne a$. Taking its adjoint gives
\begin{equation}\label{eq:normal}
 \Phi_k+\frac1{\sqrt N}\sum_{i,\ell,h}
   (u_{\ell,ik,h}C_{\ell,h,1}^*\Phi_i+
    u_{\ell,ki,h}C_{\ell,h,2}^*\Phi_i)=0,\qquad k\ne a.
\end{equation}
Left multiplication by $C_{\ell,h,1}^*$ and $C_{\ell,h,2}^*$
preserves $V$. Set $L_{\ell,h}\Phi=C_{\ell,h,1}^*\Phi$ and
$M_{\ell,h}\Phi=C_{\ell,h,2}^*\Phi$. By Lemma~\ref{lem:observable},
Theorem~\ref{thm:dual} applies with $\calD=\ker\mathcal O$
to the component of the profile in $\calD$.

\begin{proposition}[Restricted normal-space frame]\label{prop:frame}
There are $\beta<\infty$ and an event $\mathcal E_N\subset\Omega_N$
such that $\PP(\Omega_N\setminus\mathcal E_N)\le Ce^{-cN}$ and,
on $\mathcal E_N$, for every $a$ and every profile $\Phi\in V^N$
satisfying \eqref{eq:normal},
\begin{equation}\label{eq:frame}
 \frac1N\sum_{i,\ell,h,t}
       |\Tr(\Phi_i^*C_{\ell,h,t})|^2
       \ge N^{-\beta}\sum_i\norm{\Phi_i}_{\HS}^2 .
\end{equation}
\end{proposition}
\begin{proof}
Use real orthogonal projections in $V$ to write
$\Phi_i=v_i+w_i$, where $v_i\in\calD$ and
$w_i\in\calD^\perp$. Equation~\eqref{eq:normal} says
\[
 \frakD_av=-\widetilde{\frakD}_aw.
\]
Here $\frakD_a$ is the restriction to $\calD^N$, while
$\widetilde{\frakD}_a$ denotes the same expression on $V^N$.
Theorem~\ref{thm:dual} and the localized operator norm bound give,
simultaneously for all $a$,
\[
 N^{-E_*}\norm v_2
 \le\norm{\frakD_av}_2
 =\norm{\widetilde{\frakD}_aw}_2
 \le C_0\sqrt N\norm w_2.
\]
Consequently,
\[
 \norm \Phi_2^2=\norm v_2^2+\norm w_2^2
 \le (1+C_0^2N^{2E_*+1})\norm w_2^2.
\]
The restriction of $\mathcal O$ to $\calD^\perp$ is injective.
By finite dimensionality there is a constant $\lambda>0$ such that
$\norm{\mathcal O(w)}_2^2\ge\lambda\norm w^2$ there.
Using $\mathcal O(\Phi_i)=\mathcal O(w_i)$, we obtain
\[
 \frac1N\sum_i\norm{\mathcal O(\Phi_i)}_2^2
 \ge\frac{\lambda}{N}\norm w_2^2
 \ge cN^{-2E_*-2}\norm \Phi_2^2.
\]
Take $\beta=2E_*+3$ and enlarge $N_0$ to absorb $c$.
The exceptional probability is the one in Theorem~\ref{thm:dual}.
If $\calD=0$, injectivity alone proves the result with
$\mathcal E_N=\Omega_N$ and, for example, $\beta=2$.
\end{proof}

\subsection{The full determinant gradient}

Jacobi's formula identifies the normalized derivatives of
$\det T_N$ with the paired traces. Summing
Proposition~\ref{prop:frame} over the inverse rows therefore gives
a bound through the full determinant gradient. In
\cite[Sections~3.4--3.7 and~6.1]{Han}, the determinant is instead
that of a projected exposed block of fixed size. The present
summation accommodates the atoms shared by different columns.

\begin{lemma}[The determinant gradient in atom coordinates]\label{lem:gradient}
Let $F=\det T_N$ as a polynomial in the $n\mathsf d N^2$ original
real coordinates $u_{\ell,ij,h}$. It has $F(0)=1$, degree at most
$2m$ in each coordinate, and, with $R=T_N^{-1}$,
\begin{equation}\label{eq:gradient}
 \frac{\partial_{u_{\ell,ij,h}}F}{F}
 =\frac1{\sqrt N}
   \bigl(\Tr(R_{ji}A_{\ell,h})+\Tr(R_{ij}B_{\ell,h})\bigr).
\end{equation}
On $\mathcal E_N$ from Proposition~\ref{prop:frame}, after removing
a null set,
\begin{equation}\label{eq:inverse-gradient}
 \norm R\le\norm R_{\HS}
 \le N^{\beta/2}\frac{\norm{\nabla F}_2}{|F|}.
\end{equation}
\end{lemma}
\begin{proof}
The derivative matrix in the coordinate $u_{\ell,ij,h}$ is
\[
 \partial_{u_{\ell,ij,h}}T_N
 =N^{-1/2}(A_{\ell,h}\otimes E_{ij}
                   +B_{\ell,h}\otimes E_{ji}).
\]
Its rank is at most $2m$, so determinant multilinearity gives
the coordinate degree bound. Jacobi's formula then yields
\[
 \frac{\partial_{u_{\ell,ij,h}}F}{F}
 =\Tr\bigl(R\,\partial_{u_{\ell,ij,h}}T_N\bigr)
 =N^{-1/2}\bigl(\Tr(R_{ji}A_{\ell,h})
                         +\Tr(R_{ij}B_{\ell,h})\bigr).
\]
When $i=j$, the derivative matrix contains $A_{\ell,h}+B_{\ell,h}$
in that one block, so the same formula applies without an
additional diagonal factor.
In class (B) the two derivatives are
\begin{align*}
 \frac{\partial_{\Re\xi_{\ell,ij}}F}{F}
 &=N^{-1/2}\bigl(\Tr(R_{ji}A_\ell)+\Tr(R_{ij}B_\ell)\bigr),\\
 \frac{\partial_{\Im\xi_{\ell,ij}}F}{F}
 &=iN^{-1/2}\bigl(\Tr(R_{ji}A_\ell)-\Tr(R_{ij}B_\ell)\bigr).
\end{align*}
These identities hold whether or not the two atom coordinates
are independent.
Since $F(0)=1$, its zero set is Lebesgue null in that real coordinate
space. The product block law is absolutely continuous, so $T_N$ is
invertible almost surely and $\{F=0\}$ can be removed from
$\mathcal E_N$ without changing its probability.

To control the full inverse, apply \eqref{eq:frame} to $\Phi_i=\widehat R_{ai}^*$ and sum over $a$.
Equations \eqref{eq:paired-A}--\eqref{eq:paired-B} give
\[
 \sum_{a,i,\ell,h,t}|\Tr(\widehat R_{ai}C_{\ell,h,t})|^2
 =2N\frac{\norm{\nabla F}_2^2}{|F|^2},\qquad
 \sum_{a,i}\norm{\widehat R_{ai}}_{\HS}^2=2\norm R_{\HS}^2.
\]
For each $(\ell,h)$, the $t=1$ term at $(a,i)$ is $\sqrt N$
times $\partial_{u_{\ell,ia,h}}F/F$, and the $t=2$ term is
$\sqrt N$ times $\partial_{u_{\ell,ai,h}}F/F$.
Thus each derivative, including a diagonal derivative, occurs
twice in the full sum. The squared inverse norm is also doubled
by \eqref{eq:inverse-block}. Substituting the two displayed
identities into the sum of \eqref{eq:frame} gives
\[
 2\frac{\norm{\nabla F}_2^2}{|F|^2}
 \ge 2N^{-\beta}\norm R_{\HS}^2,
\]
which proves \eqref{eq:inverse-gradient} for $F=\det T_N$.
\end{proof}

\subsection{Relative polynomial anti-concentration}

The inverse bound reduces the least singular value event to a
determinant that is small relative to its gradient. We bound this
relative event directly under the atom law.

\begin{lemma}[Relative polynomial anti-concentration]\label{lem:relative}
Let $F$ be a nonzero complex-valued polynomial in real coordinates,
of degree at most $d$ in each coordinate.
\begin{enumerate}[label=(\roman*)]
\item Let $X_1,\ldots,X_k$ be independent real random variables with
densities $f_j$ satisfying $\|f_j\|_{L^\infty(\R)}\le M$.
Writing $X=(X_1,\ldots,X_k)$, we have
\begin{equation}\label{eq:relative}
 \PP\{|F(X)|\le t\norm{\nabla F(X)}_2\}
 \le 2Md^2 k^{3/2}t,\qquad t>0.
\end{equation}
\item For $k$ independent $\R^2$-valued blocks, each with joint
density bounded by $M$ and second moment at most $M_2$,
\begin{equation}\label{eq:relative-planar}
 \PP\{|F|\le t\norm{\nabla F}_2\}
 \le C_{d,M,M_2}k^{3/2}t^{1/2},\qquad 0<t\le1.
\end{equation}
Here the gradient contains all $2k$ real coordinate derivatives.
\end{enumerate}
\end{lemma}
\begin{proof}
For a nonzero univariate polynomial $q$ of degree $b\le d$,
$q'/q=\sum_{j=1}^b(x-z_j)^{-1}$ off its zeros. Thus
$|q|\le v|q'|$ forces $|x-z_j|\le bv$ for some root.
On the real line these intervals have total length at most $2d^2v$.
A nonzero constant gives an empty event, and a finite zero set
does not affect length.

For (i), on the event with $F\ne0$ some coordinate satisfies
$|F|\le t\sqrt{k}\,|\partial_j F|$. Condition on all other
coordinates and apply the length estimate with $v=t\sqrt{k}$.
Independence leaves the selected coordinate density at most $M$.
Sum over $k$ coordinates to obtain \eqref{eq:relative}.

For (ii), fix a block $(x,y)$ and condition on all the other
blocks. The resulting polynomial in $(x,y)$ is nonzero for
almost every choice of the other blocks: in the expansion in
powers of $x,y$, at least one coefficient is a nonzero polynomial
of those blocks. For a nonzero bivariate polynomial, the exceptional
$y$ for which the $x$-slice vanishes identically form a finite
set: at least one coefficient polynomial in $y$ is nonzero. All other slices obey the univariate
length estimate; a nonzero constant slice contributes nothing.
For a cutoff $R_{\mathrm{cut}}>0$, Fubini gives
\[
 \operatorname{Leb}_2\bigl\{x^2+y^2\le R_{\mathrm{cut}}^2:
      |F(x,y)|\le v|\partial_xF(x,y)|\bigr\}
 \le 4d^2R_{\mathrm{cut}}v.
\]
The same estimate holds for $\partial_yF$. Multiplying the area
by the joint density bound $M$ bounds the conditional probability.
The coordinates $x$ and $y$ may be dependent, since we
conditioned only on the other blocks.

If all $k$ blocks have norm at most $R_{\mathrm{cut}}$, the event in (ii) implies
one of the $2k$ coordinate events with $v=t\sqrt{2k}$.
The complement has probability at most $kM_2/R_{\mathrm{cut}}^2$. Therefore
\[
 \PP\{|F|\le t\norm{\nabla F}_2\}
 \le CMd^2k^{3/2}R_{\mathrm{cut}}t+kM_2R_{\mathrm{cut}}^{-2}.
\]
Take $R_{\mathrm{cut}}=t^{-1/2}$ for $0<t\le1$; then $R_{\mathrm{cut}}t=t^{1/2}$ and
$t\le t^{1/2}$, proving \eqref{eq:relative-planar}.
In both parts, all exceptional coefficient zero sets and $\{F=0\}$
are null under the product of the absolutely continuous block laws.
\end{proof}

\subsection{The least singular value cutoff}

\begin{proposition}[Irreducible pencils in both input classes]\label{prop:irreducible}
For every nonzero complex irreducible tuple,
\[
 \PP\{\smin(T_N)\le t,\ \Omega_N\}
 \le CN^Kt^\omega+Ce^{-cN}.
\]
For the Gaussian laws in Theorem~\ref{thm:pencil}, one may take
$\omega=1$ also in class (B).
\end{proposition}
\begin{proof}
On $\mathcal E_N$, the event of a small least singular value implies
$|F|\le tN^{\beta/2}\norm{\nabla F}_2$.
In class (A), apply Lemma~\ref{lem:relative}(i) with $d=2m$ and
$k=nN^2$ to get $CN^{3+\beta/2}t$.
In class (B), set $v=tN^{\beta/2}$ and use
Lemma~\ref{lem:relative}(ii) with $k=nN^2$, $d=2m$, $M_2=1$.
For $v\le1$ this gives $CN^{3+\beta/4}t^{1/2}$.
For $v>1$ the same upper bound, enlarged if necessary, is at
least one. Adding $\PP(\Omega_N\setminus\mathcal E_N)$ proves the bound.

For the Gaussian refinement in class (B), enumerate the $nN^2$
atom blocks by $j$ and write each block as
$u=Sg$, where $g$ has two independent standard real Gaussian
coordinates and the fixed $2\times2$ real matrix $S$ is invertible.
Let $\widetilde F(g)=F((Sg_j)_j)$. By the chain rule,
\[
 \norm{\nabla_g\widetilde F}_2
 \ge\smin(S)\norm{\nabla_uF}_2.
\]
Varying one coordinate of $g_j$ varies just one complex entry
by a fixed complex multiple and its conjugate entry in $X^*$
by the conjugate multiple. The matrix variation still has rank
at most $2m$. Thus $\widetilde F$ has degree at most $2m$ in
each of its $2nN^2$ independent real Gaussian coordinates.
Apply Lemma~\ref{lem:relative}(i) to $\widetilde F$ with radius
$tN^{\beta/2}/\smin(S)$. The resulting bound
$CN^{3+\beta/2}t$ therefore bounds the least singular value event
on $\mathcal E_N$.
\end{proof}

The next section combines the factor bounds. The resulting cutoff
enters the lower logarithmic tail estimate, together with a count
of intermediate small singular values.

\section{Reconstruction from irreducible factors}\label{sec:reconstruction}

The inverse estimate applies to every nonzero irreducible factor.
The passage to arbitrary coefficients follows
\cite[Section~3.1]{Han}: a composition series makes the pencil upper
triangular, and a finite inverse expansion combines its diagonal
bounds. Localization controls the off-diagonal blocks.

\begin{proof}[Proof of Theorem~\ref{thm:pencil}]
Choose a composition series of the natural $\C^m$-module for
the original complex coefficient algebra, and an orthogonal
splitting adapted to that series. Every coefficient matrix is
then block upper triangular, with
$r_{\mathrm{comp}}\le m$ diagonal factors. The nonzero factors are complex
irreducible; a factor on which every coefficient acts as zero
contributes an identity pencil.

Write $T=T_N=D+J$, with $D=\diag(D_1,\ldots,D_{r_{\mathrm{comp}}})$ block diagonal
and $J$ strictly block upper triangular in these $r_{\mathrm{comp}}$ factors. On $\Omega_N$,
$\norm J\le C\sqrt N$. Suppose that each diagonal factor has
least singular value at least $\tau_{\mathrm{fac}}\le1$, so that
$\norm{D^{-1}}\le \tau_{\mathrm{fac}}^{-1}$. The matrix $D^{-1}J$ is strictly
block upper triangular and $(D^{-1}J)^{r_{\mathrm{comp}}}=0$. Therefore
\begin{equation}\label{astyle:eq:triangular-inverse}
 T^{-1}=(I+D^{-1}J)^{-1}D^{-1}
       =\sum_{j=0}^{r_{\mathrm{comp}}-1}(-D^{-1}J)^jD^{-1}.
\end{equation}
The $j$th term has norm at most $(C\sqrt N)^j\tau_{\mathrm{fac}}^{-j-1}$.
Using $\tau_{\mathrm{fac}}\le1$ and $N\ge1$ to sum these terms gives
\[
 \norm{T^{-1}}\le C N^{(r_{\mathrm{comp}}-1)/2}\tau_{\mathrm{fac}}^{-r_{\mathrm{comp}}}.
\]
For
\[
 \tau_{\mathrm{fac}}=(2CN^{(r_{\mathrm{comp}}-1)/2}t)^{1/r_{\mathrm{comp}}}\le1,
\]
this upper bound is $(2t)^{-1}$. Hence
\[
 \{\smin(T)\le t,\ \Omega_N\}
 \subset\bigcup_{a=1}^{r_{\mathrm{comp}}}
 \{\smin(D_a)<\tau_{\mathrm{fac}},\ \Omega_N\}.
\]
Identity factors do not contribute to this union.
Proposition~\ref{prop:irreducible}, applied to the fixed nonzero
factors, bounds the right-hand side by
\[
 CN^K\tau_{\mathrm{fac}}^\omega+Ce^{-cN}
 \le CN^{K+\omega(r_{\mathrm{comp}}-1)/(2r_{\mathrm{comp}})}t^{\omega/r_{\mathrm{comp}}}+Ce^{-cN}.
\]
For Gaussian laws we use $\omega=1$. If $\tau_{\mathrm{fac}}$ exceeds one,
increasing the polynomial prefactor makes the bound at least
one. Adding $\PP(\Omega_N^c)=q_N$ gives the unconditional estimate.
If all coefficients vanish, $T=I$ and the estimate is immediate.

\paragraph{The localization error.}
Since
$N^{-1}\norm{X_\ell}_{\HS}^2=N^{-2}\sum_{ij}|\xi_{\ell,ij}|^2$,
the weak law gives $q_N\to0$. If the law has a $(2+\varepsilon)$ moment,
set $r=1+\min(\varepsilon,2)/2\in(1,2]$ and write
$Z_j=|\xi_j|^2-1$, $1\le j\le N^2$.
Let $Z_j'$ be independent copies and $\varepsilon_j$ independent
Rademacher signs. Conditional Jensen's inequality, symmetry of $Z_j-Z_j'$,
and the Rademacher second moment give
\begin{align*}
 \EE\left|\sum_jZ_j\right|^r
 &\le\EE\left|\sum_j(Z_j-Z_j')\right|^r
  =\EE\EE_\varepsilon
       \left|\sum_j\varepsilon_j(Z_j-Z_j')\right|^r\\
 &\le\EE\left(\sum_j|Z_j-Z_j'|^2\right)^{r/2}
  \le\sum_j\EE|Z_j-Z_j'|^r
  \le C_rN^2\EE\bigl||\xi|^2-1\bigr|^r.
\end{align*}
Markov's inequality at threshold $N^2$ yields
$q_N\le CN^{-2(r-1)}$.
For either Gaussian class, $|\xi|^2$ has a finite exponential
moment in a neighborhood of zero. Since $\EE|\xi|^2=1$, some
fixed $\lambda>0$ satisfies
$\log\EE e^{\lambda|\xi|^2}<2\lambda$.
The exponential Markov inequality at threshold $2N^2$, and a union over $\ell$,
give $q_N\le Ce^{-cN^2}$. Combining this with the localized
Gaussian estimate completes the proof.
\end{proof}

\section{Least singular values for Wigner tensor pencils}
\label{w:section}\label{sec:wigner}

We prove Theorem~\ref{thm:wigner-pencil} and its consequence for
polynomial shifts by matching inverse traces with determinant
derivatives in unordered-edge coordinates. We adapt the profile argument of
Section~\ref{sec:abstract} to unordered edges: an off-diagonal
atom enters both endpoint equations, so we leave only edges from
an anchor set to its complement unexposed and discard the anchor
outputs. The determinant derivative also changes both endpoint
entries. We match it to paired inverse traces, then apply relative
anti-concentration and triangular reconstruction. The resulting
cutoff will enter the Wigner lower logarithmic tail estimate.

In this section abbreviate $T_N^{\mathrm W}$ to $T_N$ and write
\[
 T_N=I_{mN}+\sum_{\ell=1}^nC_\ell\otimes H_\ell^{(N)}.
\]
The off-diagonal atoms are independent across $\ell$ and unordered
edges. For each $\ell$, they have a fixed centered law of absolute
variance one, with bounded real or joint planar density. Both types
of law may occur in the same pencil. For each $\ell$, the diagonal
atoms are independent copies of a fixed centered real law with bounded
density and finite second moment. The two coordinates of a planar atom may be dependent.

We use the following real coordinates. For $i<k$, write $\zeta_{\ell,ik}=x_{\ell,\{i,k\}}+
iy_{\ell,\{i,k\}}$, omitting $y$ for real entry laws, and set
$d_{\ell,i}=\sqrt N(H_\ell^{(N)})_{ii}$.  We use
$\varepsilon_{ik}=1$ for $i<k$ and
$\varepsilon_{ik}=-1$ for $i>k$.  Thus, for $i\ne k$,
\[
 (H_\ell^{(N)})_{ik}=N^{-1/2}
   (x_{\ell,\{i,k\}}+i\varepsilon_{ik}y_{\ell,\{i,k\}}).
\]
Let $Y_{\ell,x}$ be the real symmetric matrix with off-diagonal
entries $x_{\ell,\{i,k\}}/\sqrt N$ and diagonal
$d_{\ell,i}/\sqrt N$, and let $Y_{\ell,y}$ be the real
skew-symmetric matrix with entries
$\varepsilon_{ik}y_{\ell,\{i,k\}}/\sqrt N$ and zero diagonal.
Then $H_\ell^{(N)}=Y_{\ell,x}+iY_{\ell,y}$ and
\begin{equation}\label{w:energy}
 \|Y_{\ell,x}\|_{\HS}^2+\|Y_{\ell,y}\|_{\HS}^2
 =\|H_\ell^{(N)}\|_{\HS}^2.
\end{equation}
Throughout this section, $y$ terms are omitted for real entry laws.
In sums indexed by $h$, the value $h=x$ is always included, and
$h=y$ is included only for planar entry laws.
We write $\Omega_N^{\mathrm W}$ and $q_N^{\mathrm W}$ for the
localization event and its complementary probability in
\eqref{eq:wigner-Omega}. We use \eqref{w:energy} to bound both
the pencil norm and the contribution of small singular directions
of the profile coefficient map.

\subsection{A profile estimate with unordered edges}

The dimension surplus depends only on the coefficient maps. To
apply the uniform covering and projection arguments here, we prove
the probability estimate for whole unordered edges and include
the diagonal contribution in the approximation error.

Let $V$ be a real Hilbert space of fixed dimension $p\ge1$,
let $\calD\subset V$ be a fixed real subspace, and let $L_{\ell,x},L_{\ell,y}:V\to V$
be fixed real linear maps, with $L_{\ell,y}$ present only for
planar entry laws. Assume that no nonzero real subspace of $\calD$
is invariant under all the available maps. Define, for
$a\in[N]$ and $v\in\calD^N$,
\begin{equation}\label{w:profile}
 \begin{split}
 (\frakD_a^{\mathrm W}v)_k={}&v_k+
 \frac1{\sqrt N}\sum_{i\ne k,\ell}
 \bigl(x_{\ell,\{i,k\}}L_{\ell,x}v_i+
       \varepsilon_{ik}y_{\ell,\{i,k\}}L_{\ell,y}v_i\bigr)\\
 &+\frac1{\sqrt N}\sum_\ell d_{\ell,k}L_{\ell,x}v_k,
 \qquad k\ne a.
 \end{split}
\end{equation}
The same formula on $V^N$ defines
$\widetilde{\frakD}_a^{\mathrm W}$.  Equation~\eqref{w:energy}
implies, on $\Omega_N^{\mathrm W}$, that
\begin{equation}\label{w:profile-norm}
 \|\widetilde{\frakD}_a^{\mathrm W}\|\le C_0\sqrt N
 \quad\hbox{simultaneously for all }a.
\end{equation}
The same localization gives a bound on the coordinate squares:
\[
 \sum_\ell\left(2\sum_{i<k}(x_{\ell,\{i,k\}}^2+
        y_{\ell,\{i,k\}}^2)+\sum_i d_{\ell,i}^2\right)\le2nN^2.
\]
This also bounds the inactive part of the profile map, including
the diagonal terms. Constants below depend on the fixed maps and laws.

\begin{lemma}[Unordered-edge profile estimate]\label{w:dual}
If $\calD\ne0$, there are $E_*<\infty$ and $N_0$ such that, for $N\ge N_0$,
\[
 \PP\left(\Omega_N^{\mathrm W}\cap
   \{\exists a:\smin(\frakD_a^{\mathrm W})\le N^{-E_*}\}\right)
 \le Ce^{-cN}.
\]
The conclusion uses no moment above two and no independence within
a planar off-diagonal atom.
\end{lemma}
\begin{proof}
For a nonzero $E\subset\calD$, let $r$ be the number of available
maps $L_{\ell,x},L_{\ell,y}$ and let $\Psi_E:E^r\to V$ sum them. Define its active
image $Y_\sigma(E)$ and the right singular space
$C_{\tau,\sigma}(E)$ of $P_{Y_\sigma(E)^\perp}|_E$ as in
Section~\ref{sec:abstract}. Lemma~\ref{lem:surplus}, applied
to this list of maps, gives fixed $\sigma_0,\tau>0$ with
\begin{equation}\label{w:surplus}
 y:=\dim Y_\sigma(E)\ge
 c_E+1,\qquad c_E:=\dim C_{\tau,\sigma}(E),
 \quad 0<\sigma\le\sigma_0.
\end{equation}
All dimensions here are real.

\emph{Conditioning on whole edges.}
Fix a deterministic
$v\in E^N$ with
$\|v\|_2^2\le C_1$ and $\Gamma_v\succeq\gamma P_E$,
where $d=\dim E$.  Lemma~\ref{lem:anchors} selects, before any
randomness is exposed, an anchor set $I$ with $|I|=d$ and
\[
 \Gamma_{v,I}:=\sum_{i\in I}v_i\otimes_{\R}v_i
 \succeq c_{d,C_1}\gamma^dN^{-d}P_E.
\]
Let $J=[N]\setminus(I\cup\{a\})$ and expose every atom except
the whole edge blocks on $\{i,k\}$ with $i\in I$ and $k\in J$.
In particular, all diagonal atoms are exposed. The remaining random
part of output $k\in J$ is
\[
 Z_k=N^{-1/2}\sum_{i\in I,\ell}
 \bigl(x_{\ell,\{i,k\}}L_{\ell,x}v_i+
       \varepsilon_{ik}y_{\ell,\{i,k\}}L_{\ell,y}v_i\bigr).
\]
The second occurrence of every unexposed block is in an anchor
output, which has been discarded.  The remaining outputs
therefore use disjoint independent blocks.  All other terms of
their equations are fixed translations.

For $Y=Y_\sigma(E)$, the real linear map from the unexposed
coordinates to $P_YZ_k$ has Euclidean Gram operator
\[
 \frac1N P_Y\sum_\ell
 \bigl(L_{\ell,x}\Gamma_{v,I}L_{\ell,x}^*
       +L_{\ell,y}\Gamma_{v,I}L_{\ell,y}^*\bigr)P_Y
 \succeq c\sigma^2\gamma^dN^{-d-1}I_Y.
\]
Each $y$ coefficient is multiplied by $\varepsilon_{ik}$,
whose square is one, so its outer product is the same at either
endpoint. The lower bound follows from the anchor estimate
and $\Psi_E\Psi_E^*\succeq\sigma^2I_Y$ on $Y$. This calculation
concerns the deterministic coordinate map, so dependence between
the coordinates of a planar block does not affect it.

Fix $\kappa=1/(8p)$. Since only finitely many laws occur, we
can choose $R_0$ so large that, for every anchor set of size at
most $p$, all its $n|I|$ edge blocks in one output have norm at
most $R_0$ with probability at least $1-\kappa/4$. This choice
uses only that each law is a probability measure. For a fixed
anchor set, the good-output indicators involve disjoint edges and
are independent. Chernoff's inequality, followed by a union over
at most $pN^p$ anchor sets, gives an event
$\mathcal A_N^{\mathrm W}$ with
$\PP((\mathcal A_N^{\mathrm W})^c)\le Ce^{-cN}$, on which
at most $\kappa N$ outputs fail the cutoff for any anchor set.

Condition on the pattern of outputs satisfying the cutoff. In each
such output the unexposed edge blocks are independent copies of their
original laws conditioned to lie in the radius-$R_0$ ball. Each
planar block is conditioned as a whole. At most $np$ blocks occur,
each of real dimension one or two, so their product density and
support radius are bounded by fixed constants. The preceding Gram
bound and Lemma~\ref{lem:projection} then give the one-output
small-ball factor
\[
 C\sigma^{-y}\gamma^{-dy/2}N^{(d+1)y/2}\delta^y.
\]
On $\mathcal A_N^{\mathrm W}$ the chosen anchor set has at least
$b_N(d)=\lfloor(1-\kappa)N\rfloor-d-1$ good outputs in $J$.
Replacing the global cutoff event by this assertion for the chosen
$I$, we multiply the one-output factors and average over the patterns. Their
probabilities sum to at most one, so we obtain
\begin{equation}\label{w:pointwise}
 \PP\{\|\frakD_a^{\mathrm W}v\|_2\le\delta,
                 \mathcal A_N^{\mathrm W}\}
 \le\min\left\{1,
 \left(C\sigma^{-y}\gamma^{-dy/2}
       N^{(d+1)y/2}\delta^y\right)^{b_N(d)}\right\},
\end{equation}
Integrating over the first exposure proves the estimate for the
intersection with $\mathcal A_N^{\mathrm W}$ without conditioning
on that event; no restriction by $\Omega_N^{\mathrm W}$ was used.

\emph{Uniformity over profiles.}
The fixed-profile bound has the same exponent and scale factors as
in the independent model. To apply the subspace and tube nets of
Lemma~\ref{lem:regular}, we bound the remaining contribution in
the inactive singular directions.
At a net center $E$ in the family with a
spectral gap, the random sum in output $k$ is $\Psi_Et_k$, where
\[
 \begin{split}
 t_{\ell,x,k}&=N^{-1/2}\left(
       \sum_{i\ne k}x_{\ell,\{i,k\}}v_i+d_{\ell,k}v_k\right),\\
 t_{\ell,y,k}&=N^{-1/2}\sum_{i\ne k}
             \varepsilon_{ik}y_{\ell,\{i,k\}}v_i.
 \end{split}
\]
Cauchy--Schwarz and the preceding bound on the coordinate squares,
including the diagonal terms, give
\[
 \begin{split}
 \norm{t_k}^2
 &\le\frac{\norm v_2^2}{N}\sum_\ell
       \left(\sum_{i\ne k}(x_{\ell,\{i,k\}}^2+
                           y_{\ell,\{i,k\}}^2)+d_{\ell,k}^2\right)\\
 &\le2nN\norm v_2^2.
 \end{split}
\]
Thus projection of the random sum onto the orthogonal complement
of the active image leaves a contribution of at most $C\sqrt N\eta$. When this is smaller
than the small-image radius, the profile coordinates are confined
to the same tubes as before. The operator bound
\eqref{w:profile-norm} controls rounding errors, while
\eqref{w:surplus} gives the excess of active output dimensions
over the dimensions of the profile net. The same Grassmannian
and tube construction therefore applies, including its Gram
stability estimate. Summing the integrated probabilities
\eqref{w:pointwise} over this net leaves the dimension surplus
$yb_N(d)-c_E(N-1)-d\ge7N/8-O_p(1)$, yielding
\begin{equation}\label{w:uniform-net}
 C_*^N N^{K_0N}\gamma^{-p^2N}\sigma^{-pN}
                   \delta^{N/2-L_0}
\end{equation}
for the probability of a small image among regular profiles in
$\Omega_N^{\mathrm W}\cap\mathcal A_N^{\mathrm W}$, under the
same conditions $0<\delta\le c_*\gamma$ and
$C_*\sqrt N\eta\le\delta$. The constants depend on the fixed
maps and entry laws.

\emph{Spectral scales and projection of profiles.}
The remaining reductions use only the uniform probability bound
and the localized operator norm. Applying the spectral scales of Lemma~\ref{lem:all-regular}
to \eqref{w:uniform-net} and the surplus gives
a cutoff $N^{-e(D)}$ for
$\Gamma_v\succeq N^{-D}P_E$, with failure probability
$pe^{-2N\log N}$. We then project away the eigenvectors corresponding
to small Gram eigenvalues
as in the proof of Theorem~\ref{thm:dual}, using the
$C_0\sqrt N$ operator bound. With the
thresholds $D_1=1$,
$D_{d+1}\ge\max\{D_d+1,2e(D_d)+4\}$, the projection error
in the image is at most $\tfrac12N^{-e(D_d)}$ at the stopping
dimension. Hence $E_*>\max_d e(D_d)+1$ suffices. Taking a union
over the finitely many dimensions and the $N$ choices of $a$, and
adding $\PP((\mathcal A_N^{\mathrm W})^c)$ proves the lemma.
\end{proof}

\subsection{Mixed real and planar relative anti-concentration}

The inverse estimate uses the real gradient over independent atom
blocks. Off-diagonal atoms may be planar, while diagonal atoms
are real, so the relative estimate must allow both types. For a
planar block, we bound the area of the small-value set before
applying its joint density, without conditioning one coordinate
on the other.

\begin{lemma}\label{w:mixed-relative}
Let $F$ be a nonzero complex-valued polynomial of degree at most
$d$ in each real coordinate.  Its independent atom blocks comprise
$k_1$ real variables with densities bounded by $M$, and $k_2$
planar variables with joint densities bounded by $M$ and second
moments at most $M_2$.  Write $K_{\mathrm{coord}}=k_1+2k_2\ge1$.  Then
\[
 \PP\{|F|\le t\|\nabla F\|_2\}
 \le C_{d,M,M_2}K_{\mathrm{coord}}^{3/2}t^{1/2},\qquad0<t\le1.
\]
If $k_2=0$, one can replace $t^{1/2}$ by $t$.
\end{lemma}
\begin{proof}
Outside the null zero set of $F$, the event in question implies,
for some real coordinate $j$,
\[
 |F|\le t\sqrt{K_{\mathrm{coord}}}\,|\partial_jF|.
\]
For a real atom, fix the other whole blocks and use the
root-factorization estimate from Lemma~\ref{lem:relative}.
The conditional probability is at most $2Md^2t\sqrt{K_{\mathrm{coord}}}$.
For a coordinate in a planar atom, fix only the other blocks.
Within the radius-$R_{\mathrm{cut}}$ disk, the one-variable root estimate bounds
every slice by length $2d^2t\sqrt{K_{\mathrm{coord}}}$. Fubini therefore bounds
the admissible area by $4d^2R_{\mathrm{cut}}t\sqrt{K_{\mathrm{coord}}}$. Multiplying by the joint
density bound gives the conditional probability estimate without
assuming independence of the two coordinates. Summing over the coordinates and using the
second-moment bound outside the planar disks gives
\[
 \PP\{|F|\le t\|\nabla F\|_2\}
 \le C Md^2 K_{\mathrm{coord}}^{3/2}(t+R_{\mathrm{cut}}t)+k_2M_2R_{\mathrm{cut}}^{-2}.
\]
Taking $R_{\mathrm{cut}}=t^{-1/2}$ proves the square-root bound. The exceptional
coefficient and slice zero sets are null under the product of the
absolutely continuous block laws, as in Lemma~\ref{lem:relative}.
When there are no planar atoms, the real-variable estimate alone
gives the stated linear bound.
\end{proof}

\subsection{The doubled algebra and the determinant gradient}

We choose the doubled generators so that their traces encode the
symmetric and skew-symmetric variations of an unordered edge.
The profile estimate gives a trace frame on the adjoint algebra.
We identify each trace with an edge or diagonal derivative and sum
the frame inequalities over the rows.

Let $(C_\ell)$ be a nonzero complex irreducible coefficient family.
Since an edge derivative contains the sum or difference of the
traces at its endpoints, define
\begin{equation}\label{w:doubled-coefficients}
 D_{\ell,x}=\diag(C_\ell,C_\ell^{\mathsf T}),\qquad
 D_{\ell,y}=\diag(iC_\ell,-iC_\ell^{\mathsf T}),\qquad
 \calB=\alg_\C(I,D_{\ell,x},D_{\ell,y}),
\end{equation}
again omitting $y$ for real entry laws. The associated
pencil satisfies
\begin{equation}\label{w:double-pencil}
 \widehat T_N=I+\sum_\ell
       (D_{\ell,x}\otimes Y_{\ell,x}
        +D_{\ell,y}\otimes Y_{\ell,y})
 \simeq T_N\oplus T_N^{\mathsf T}.
\end{equation}
The upper block is $T_N$; the lower one is $T_N^{\mathsf T}$
because $Y_{\ell,x}^{\mathsf T}=Y_{\ell,x}$ and
$Y_{\ell,y}^{\mathsf T}=-Y_{\ell,y}$.

The two coordinate representations of $\calB$ are generated by
$C_\ell$ and $C_\ell^{\mathsf T}$, up to scalar multiples.
Both are irreducible, the second by the annihilator argument in
Lemma~\ref{lem:algebra}. The same lemma gives a nondegenerate
trace form and zero common kernel for these doubled generators.
On the real Hilbert space $V=\calB^*$ define
\[
 \mathcal O(\Phi)=\bigl(\Tr(\Phi^*D_{\ell,h})\bigr)_{\ell,h},
 \qquad\calD=\ker\mathcal O.
\]
Left multiplication preserves $V$, since
$D_{\ell,h}^*b^*=(bD_{\ell,h})^*$. The right-ideal argument
of Lemma~\ref{lem:observable} therefore applies and gives \eqref{eq:no-invariant}
for this $\calD$.

If $R=T_N^{-1}$, the spatial blocks of the doubled inverse are
\[
 \widehat R_{ij}=\diag(R_{ij},R_{ji}^{\mathsf T})\in\calB.
\]
Membership follows from the Cayley--Hamilton argument of
Lemma~\ref{lem:inverse-blocks}, which requires only that the
spatial blocks of the doubled pencil belong to the unital algebra.
Taking adjoints in
$\widehat R\widehat T_N=I$ shows that
$\Phi_i=\widehat R_{ai}^*$ solves \eqref{w:profile} with zero
outputs for $k\ne a$, with $L_{\ell,h}\Phi=D_{\ell,h}^*\Phi$.
For this profile, projection onto the kernel of the trace map and
its orthogonal complement gives the same decomposition as in
Proposition~\ref{prop:frame}. The part in the kernel is controlled
by Lemma~\ref{w:dual}; the contribution of its complement to the
equation is bounded by \eqref{w:profile-norm}. Injectivity of
the trace map on that fixed complement then gives a finite $\beta$ and an event
$\mathcal E_N^{\mathrm W}\subset\Omega_N^{\mathrm W}$,
with $\PP(\Omega_N^{\mathrm W}\setminus\mathcal E_N^{\mathrm W})
\le Ce^{-cN}$, on which every inverse-row profile satisfies
\begin{equation}\label{w:frame}
 \frac1N\sum_{i,\ell,h}|\Tr(\Phi_i^*D_{\ell,h})|^2
 \ge N^{-\beta}\sum_i\|\Phi_i\|_{\HS}^2.
\end{equation}
If $\calD=0$, finite-dimensional injectivity alone gives this
estimate, with no use of the dual lemma.

The trace observations control each doubled inverse row, so summing
their squared moduli controls the full Hilbert--Schmidt inverse norm. We express
them through normalized derivatives of the original determinant;
the doubled pencil serves only to pair the endpoint traces.
Let $F=\det T_N$, viewed as a polynomial in the real coordinates
of the independent atom blocks. The gradient contains one derivative
for each available coordinate of each unordered edge and each diagonal
atom. An $x$ coordinate changes the two endpoint
entries with the same sign, a $y$ coordinate changes them with
opposite imaginary signs, and a diagonal coordinate changes only
one entry. Jacobi's formula consequently gives, for $i<j$,
\begin{align}
 \frac{\partial_{x_{\ell,\{i,j\}}}F}{F}
 &=N^{-1/2}\bigl[\Tr(R_{ji}C_\ell)+\Tr(R_{ij}C_\ell)\bigr],
       \label{w:gradient-x}\\
 \frac{\partial_{y_{\ell,\{i,j\}}}F}{F}
 &=iN^{-1/2}\bigl[\Tr(R_{ji}C_\ell)-\Tr(R_{ij}C_\ell)\bigr],
       \label{w:gradient-y}\\
 \frac{\partial_{d_{\ell,i}}F}{F}
 &=N^{-1/2}\Tr(R_{ii}C_\ell).\label{w:gradient-diag}
\end{align}
Each coordinate changes $T_N$ by rank at most $2m$, so the
degree of $F$ in that coordinate is at most $2m$. Since
$F(0)=1$, its zero set is null under the product of the continuous
atom laws.

To sum these observations, write $A=\Tr(R_{ai}C_\ell)$ and
$B=\Tr(R_{ia}C_\ell)$. Direct multiplication of the doubled
blocks gives
\[
 \Tr(\widehat R_{ai}D_{\ell,x})=A+B,\qquad
 \Tr(\widehat R_{ai}D_{\ell,y})=i(A-B).
\]
For $a\ne i$, these are $\sqrt N$ times the corresponding
unordered-edge derivatives divided by $F$, up to the endpoint
sign of the $y$ term. Each edge occurs twice in the sum over
$(a,i)$. For $a=i$, the $x$ observation is
$2\Tr(R_{ii}C_\ell)$ and the $y$ observation vanishes, whereas
the diagonal derivative has only one copy of
$\Tr(R_{ii}C_\ell)$. Write $\nabla_{\mathrm{off}}$ and
$\nabla_{\mathrm{diag}}$ for the derivatives in the unordered-edge
and diagonal coordinates, respectively. Thus the full sum $S$ of
squared observations satisfies
\begin{equation}\label{w:gradient-sum}
 S=2N\|\nabla_{\mathrm{off}}F/F\|_2^2
       +4N\|\nabla_{\mathrm{diag}}F/F\|_2^2
 \le4N\|\nabla F/F\|_2^2,
 \qquad
 \sum_{a,i}\|\widehat R_{ai}\|_{\HS}^2=2\|R\|_{\HS}^2.
\end{equation}
Sum \eqref{w:frame} over $a$ and use
\eqref{w:gradient-sum} to obtain, on
$\mathcal E_N^{\mathrm W}$ outside the null singular set,
\begin{equation}\label{w:inverse-gradient}
 \|T_N^{-1}\|\le\|T_N^{-1}\|_{\HS}
 \le\sqrt2 N^{\beta/2}\|\nabla F\|_2/|F|.
\end{equation}

\subsection{Completion of the Wigner least singular value theorem}

Relative anti-concentration now gives the tail on each irreducible
factor. For arbitrary coefficients, localization bounds the
off-diagonal blocks and the finite inverse expansion combines the
factor estimates. The localization error and Gaussian improvement
are then estimated as in the independent model.

\begin{proof}[Proof of Theorem~\ref{thm:wigner-pencil}]
For a nonzero irreducible family, \eqref{w:inverse-gradient}
implies that $\smin(T_N)\le t$ on $\mathcal E_N^{\mathrm W}$
only if
\[
 |F|\le\sqrt2\,tN^{\beta/2}\norm{\nabla F}_2.
\]
There are $O(N^2)$ independent atom blocks. If all
off-diagonal entry laws are real, Lemma~\ref{lem:relative}(i)
bounds this probability by $CN^Kt$; otherwise
Lemma~\ref{w:mixed-relative} bounds it by $CN^Kt^{1/2}$.
For relative radii above one the same conclusion follows from the
trivial probability bound after increasing the polynomial
prefactor. Thus, with
$\omega=1$ when all entry laws are real and $\omega=1/2$ otherwise,
\[
 \PP\{\smin(T_N)\le t,\Omega_N^{\mathrm W}\}
 \le CN^Kt^\omega+Ce^{-cN}.
\]

For arbitrary coefficients, choose a composition series of the
natural $\C^m$-module for their unital algebra and a fixed
orthogonal splitting. The coefficients are block upper triangular
with $r_{\mathrm{comp}}\le m$ simple diagonal factors. Factors on which all
coefficients vanish give identity pencils; the others satisfy the preceding bound on
$\Omega_N^{\mathrm W}$.

On $\Omega_N^{\mathrm W}$ the strictly upper triangular part
has norm $O(\sqrt N)$. Thus the inverse expansion
\eqref{astyle:eq:triangular-inverse} in Section~\ref{sec:reconstruction}
gives, when all diagonal least singular values are at least $\tau_{\mathrm{fac}}\le1$,
\[
 \norm{T_N^{-1}}\le CN^{(r_{\mathrm{comp}}-1)/2}\tau_{\mathrm{fac}}^{-r_{\mathrm{comp}}},\qquad
 \tau_{\mathrm{fac}}=(2CN^{(r_{\mathrm{comp}}-1)/2}t)^{1/r_{\mathrm{comp}}}.
\]
For this radius the bound is $(2t)^{-1}$, so a finite union over
the diagonal factors gives
$CN^{K'}t^{\omega/r_{\mathrm{comp}}}+Ce^{-cN}$ on $\Omega_N^{\mathrm W}$.
For larger radii the trivial bound suffices; if all coefficients
vanish, $T_N=I$. Adding $q_N^{\mathrm W}$ once proves the tail.

Finally,
\[
 N^{-1}\|H_\ell^{(N)}\|_{\HS}^2
 =N^{-2}\left(2\sum_{i<j}|\zeta_{\ell,ij}|^2
                    +\sum_i d_{\ell,i}^2\right).
\]
As in the proof of Theorem~\ref{thm:pencil}, the weak law applies to the
unordered off-diagonal sum; the diagonal term is negligible by
Markov, since its expectation is $O(N^{-1})$. Thus
$q_N^{\mathrm W}\to0$. With $(2+\varepsilon)$ moments, apply
the estimate for centered sums from that proof separately to these two sums
at threshold $cN^2$. Their normalized mean is $1+O(N^{-1})$,
below the localization threshold, so
$q_N^{\mathrm W}\le CN^{-\min(\varepsilon,2)}$.

For Gaussian laws, the whitening and chain-rule argument of
Proposition~\ref{prop:irreducible} applies block by block, including
the real diagonal atoms. A whitened edge coordinate is still
supported on its two endpoints, so its variation has rank at most
$2m$; diagonal variations have rank at most $m$. The resulting
independent real coordinates give $\omega=1$. The exponential Markov
argument in the proof of Theorem~\ref{thm:pencil} gives
$q_N^{\mathrm W}\le Ce^{-cN^2}$. The inverse expansion gives
radius exponent $1/r_{\mathrm{comp}}$ and total additive error $Ce^{-cN}$.
\end{proof}

For polynomials in the Hermitian inputs, the word-chain
linearization gives the following consequence.

\begin{corollary}[Polynomial shifts in the Wigner model]
\label{w:polynomial-lsv}
For a fixed scalar-coefficient polynomial $P$ in the Hermitian
variables and every fixed $z\ne P(0)$, finite second moments and
the stated bounded densities imply
\[
 \PP\{\smin(P(H^{(N)})-zI)\le t\}
 \le C_zN^{a_z}t^{b_z}+C_ze^{-c_zN}+q_N^{\mathrm W}
\]
for some fixed positive constants.  In particular a polynomial
least-singular-value cutoff holds with probability tending to one.
\end{corollary}
\begin{proof}
The constant case is immediate. Otherwise Lemma~\ref{lem:linearization},
with $Z_\ell=H_\ell^{(N)}=Z_\ell^*$, gives the word-chain pencil
\[
 \mathcal L_z(H)=K_z\otimes I+
       \sum_\ell(U_\ell+V_\ell)\otimes H_\ell^{(N)},
 \qquad \det K_z=P(0)-z\ne0.
\]
Multiplication by $K_z^{-1}\otimes I$ gives a Wigner pencil with
identity constant term. The word-chain elimination identity expresses
$(P(H)-zI)^{-1}$ as the upper-left block of
$\mathcal L_z(H)^{-1}$, and hence
\[
 \|(P(H)-zI)^{-1}\|\le\|\mathcal L_z(H)^{-1}\|
 \le\|[(K_z^{-1}\otimes I)\mathcal L_z(H)]^{-1}\|
          \|K_z^{-1}\|.
\]
The tail estimate for the normalized pencil proves the assertion.
The comparison uses only the fixed matrix $K_z^{-1}$ and requires
no norm estimate for the elimination factors.
\end{proof}

\section{Logarithmic potentials and Brown measure convergence}\label{sec:spectral}

We prove Theorem~\ref{thm:brown} by comparing normalized logarithmic
determinants, then take distributional Laplacians to identify the
eigenvalue measure.
Fix a nonconstant polynomial $P$.
For a fixed shift $z$, let $\nu_N^z$ be the squared singular-value law of
$P(X,X^*)-zI$, and let $\nu^z$ be the law of
$(P(c,c^*)-z)(P(c,c^*)-z)^*$. In the iid argument, $P(Z)$ abbreviates
$P(Z,Z^*)$; arguments are omitted when the input tuple is fixed.
Weak convergence alone does not justify
the logarithmic integral. We need both logarithmic tails and convergence
after planar integration to identify the Brown measure. We prove:
\begin{enumerate}[label=(\roman*)]
\item $\nu_N^z$ converges weakly in probability to $\nu^z$ under a
finite second moment.
\item For $z\ne P(0)$, the logarithmic integrals converge in probability.
The finite $(2+\varepsilon)$ moment gives a count down to a polynomial
scale; Corollary~\ref{cor:poly} controls the interval next to zero.
Together they bound the negative logarithmic tail.
\item The logarithmic potentials converge in expected local spatial
$L^1$. Taking their distributional Laplacians then gives Brown measure
convergence.
\end{enumerate}
We extend the Gaussian and non-Gaussian Hermitization arguments of
\cite[Sections~5 and~6.2]{Han}, based on the logarithmic-potential method
of \cite{TVK,BC}. Rank truncation and Taylor remainder estimates allow
the stated weaker moment assumptions for real and correlated planar atoms.

The least singular value estimates in
Sections~\ref{sec:abstract}--\ref{sec:reconstruction} control the interval
closest to zero. The count above that cutoff starts from two results
for the free and Gaussian laws.
A nonconstant self-adjoint scalar polynomial in
free semicircular variables has a non-atomic law and an algebraic Cauchy
transform with local Puiseux expansions
\cite[Theorem 1.1, parts (2), (5), (6), and Corollary 1.2]{SS}.
For a fixed self-adjoint polynomial in independent normalized GUE
matrices, the expected Cauchy transform at $i\eta$ differs from its free
limit by $CN^{-2}\eta^{-b_1}$ for $N^{-b_0}\le\eta\le1$, for some
$b_0,b_1>0$ \cite[Lemma 5.5.4]{AGZ}. We apply the latter result only to
GUE matrices, writing a standard complex Ginibre matrix as
$(H_1+iH_2)/\sqrt2$. We compare the present entry laws with these
Gaussian laws below.

\subsection{Rank truncation and weak singular-value convergence}
\label{subsec:weak-log}

A rank bound will control the positive logarithmic tail under second
moments. Applied to the elimination matrices, it will also transfer the
small-value count to the polynomial. The same truncation allows
approximation by bounded atoms.
For a matrix $\mathsf M$, let $\Pi_T\mathsf M$
replace each singular value $s$ by $\min(s,T)$. This is the
Hilbert--Schmidt projection onto the convex operator-norm ball. Hence
\[
 \norm{\Pi_T\mathsf M-\Pi_T\mathsf M'}_{\HS}\le\norm{\mathsf M-\mathsf M'}_{\HS},\qquad
 \rank(\mathsf M-\Pi_T\mathsf M)\le T^{-2}\norm{\mathsf M}_{\HS}^2.
\]
For Hermitian $\mathsf M$, this clips the eigenvalues to $[-T,T]$.
We apply the truncation separately to each input matrix.

\begin{lemma}[Upper singular-value count from second moments]
\label{sstyle:lem:rank-cap}
Suppose $Z=(Z_1,\ldots,Z_n)$ satisfies
$\EE\norm{Z_\ell}_{\HS}^2\le C N$. If $Q$ is a fixed square matrix-valued
$*$-polynomial of degree at most $D\ge1$, whose coefficient matrices
have fixed dimensions, then
\begin{equation}\label{eq:upper-count}
 \EE\#\{j:s_j(Q(Z))>M\}\le C_Q NM^{-2/D},\qquad M\ge1.
\end{equation}
The constants may depend on the fixed coefficient dimensions and the
second-moment bound.
\end{lemma}
\begin{proof}
For each word, replacing its factors one at a time gives
\[
 \rank(Q(Z)-Q(Z'))\le C_Q\sum_\ell\rank(Z_\ell-Z'_\ell).
\]
Indeed, telescoping a word gives summands with a factor
$Z_\ell-Z'_\ell$ or its adjoint, and multiplication cannot increase
rank. With $Z'=\Pi_TZ$, it follows that
\[
 \EE\rank(Q(Z)-Q(\Pi_TZ))\le C_QNT^{-2},\qquad
 \norm{Q(\Pi_TZ)}\le C_Q\max(1,T)^D.
\]
For all sufficiently large $M$, choose $T=c_QM^{1/D}$ so that the
last norm is at most $M$. Rank interlacing then bounds the number of
singular values of $Q(Z)$ above $M$ by the rank difference. This proves
\eqref{eq:upper-count} for large $M$; increasing the constant gives the
remaining values $M\ge1$. The estimate also holds for constant polynomials.
\end{proof}

For assertion (i), suppose first that the atoms are bounded, centered,
and have absolute variance one. Expand the
normalized trace of a fixed $*$-word of length
$k$ as a sum over closed index walks. Each term has normalization
$N^{-1-k/2}$. If it involves $e$ distinct atoms $(\ell,a,b)$ and
$v$ spatial indices, centering and connectedness imply
\[
 e\le k/2,\qquad v\le e+1\le1+k/2.
\]
There are finitely many equality patterns for fixed $k$, and boundedness
controls the atom moments in their contributions.
Therefore only patterns with $v=1+k/2$ can have a nonzero limit.

Equality forces $e=k/2$ and the spatial graph to be a tree.
Loops, parallel edges with different matrix labels, or counting the
reverse entry as a separate atom would reduce the number of vertices. Each edge of
the tree is traversed exactly twice, once in each direction. For the
two factors to use the same directed atom, they must have opposite
star status, so their expectation is $\EE(\xi\bar\xi)=1$.
The number of distinct vertex labelings is
$N(N-1)\cdots(N-v+1)$, and its normalized contribution tends to one.
Leaf removal gives precisely the circular noncrossing pairings
with matching matrix indices and opposite star status.
A pairing involving $\EE\xi^2=\rho$
or its conjugate would have equal star status and hence cannot occur
in a leading tree pattern. Every remaining nonzero pattern is
$O(N^{-1/2})$.

In the covariance of two normalized traces, patterns with disjoint
atom sets cancel. A shared atom joins the two underlying graphs, and a
surviving pattern of total length $k$ still has
$v\le e+1\le1+k/2$. Its normalization is now $N^{-2-k/2}$, so its
contribution is $O(N^{-1})$. The variance therefore tends to zero,
and the joint $*$-moments converge in probability. Applying
this convergence to powers of $(P-z)(P-z)^*$ gives convergence of every moment of
$\nu_N^z$ in probability. To pass to weak convergence, take an arbitrary
subsequence and a further subsequence along which all these moments
converge almost surely. The next higher moments give tightness and
uniform integrability of each lower moment. Every subsequential weak
limit therefore has the moments of $\nu^z$. Since $\nu^z$ has compact
support, moment determinacy gives weak convergence in probability.

For atoms with finite second moments, define the bounded approximants by
\[
 a_L=\EE[\xi\one_{\{|\xi|\le L\}}],\qquad
 b_L^2=\EE|\xi\one_{\{|\xi|\le L\}}-a_L|^2,\qquad
 \xi_L=b_L^{-1}(\xi\one_{\{|\xi|\le L\}}-a_L).
\]
For all sufficiently large $L$, $b_L>0$, and
\[
 a_L\longrightarrow0,\qquad b_L\longrightarrow1,\qquad
 \delta_L^2:=\EE|\xi-\xi_L|^2\longrightarrow0.
\]
The bounded law $\xi_L$ is centered and has absolute variance one.
Its pseudovariance may differ from that of $\xi$; the preceding
moment argument applies in either case.

Couple $X$ and $Y_L$ entrywise through $\xi$ and $\xi_L$. This changes
most entries, so their difference need not have small rank. Insert
operator-norm truncations between the two polynomials and compare
\[
 P(X)-zI,\quad P(\Pi_TX)-zI,\quad
 P(\Pi_TY_L)-zI,\quad P(Y_L)-zI.
\]
The two outer differences have expected rank at most $C_PNT^{-2}$.
For the middle difference, telescope each word and use the norm bounds
on the truncated inputs and nonexpansiveness of $\Pi_T$ to obtain
\[
 \EE N^{-1}\norm{P(\Pi_TX)-P(\Pi_TY_L)}_{\HS}^2
 \le C_{P,z}(T)\sum_\ell
       \EE N^{-1}\norm{X_\ell-Y_{\ell,L}}_{\HS}^2
 \le C_{P,z}(T)\delta_L^2.
\]
Both middle outputs have norm at most a fixed $C_{P,z}(T)$.
Combining Mirsky's inequality and Cauchy--Schwarz with
$|s^2-t^2|\le2C_{P,z}(T)|s-t|$ bounds the difference of their squared
singular-value laws. Write $\nu_{Z,N}^z$ for the squared singular-value
law with input $Z$. Then, for $f\in C_c^1(\R)$,
\begin{equation}\label{sstyle:eq:truncation-test}
 \EE\left|\int f\dd\nu_{X,N}^z-\int f\dd\nu_{Y_L,N}^z\right|
 \le C_fT^{-2}+C_{f,P,z}(T)\delta_L.
\end{equation}
For the outer comparisons, integrate $f'$ against the distribution
functions and use rank interlacing. For the middle comparison, use
\[
 \frac1N\sum_j|s_j^2-t_j^2|
 \le2C_{P,z}(T)
       \left(\frac1N\sum_j|s_j-t_j|^2\right)^{1/2}
 \le2C_{P,z}(T)N^{-1/2}\norm{\mathsf A-\mathsf B}_{\HS},
\]
where $\mathsf A,\mathsf B$ are the middle outputs. This proves
\eqref{sstyle:eq:truncation-test}.
Choose first $T$ large and then $L$ large, and apply the bounded-law
limit with this fixed $L$ as $N\to\infty$. The upper count
\eqref{eq:upper-count} gives tightness, so these compactly supported
test functions suffice, proving assertion (i). To integrate the negative
logarithm, we next need a quantitative bound on mass near zero.

\subsection{Small singular values for the free and Gaussian laws}\label{subsec:gaussian-reference}

For assertion (ii), assume now that $\EE|\xi|^{2+\varepsilon}<\infty$
for some $\varepsilon>0$. We need a power bound for the expected
small-singular-value count down to a polynomial scale. We first prove it
for the free law and the Gaussian polynomial, then for its linearization.
The pencil and elimination matrices are evaluated as matrix-valued
$*$-polynomials:
$\mathcal L_z(G)$ uses $G^*$, and on real inputs it agrees with the
transpose pencil.

Write each circular variable as
$c_\ell=(s_{\ell,1}+is_{\ell,2})/\sqrt2$. The scalar self-adjoint
polynomial $Q_z=(P-z)(P-z)^*$ in these semicircular generators is
nonconstant. Indeed, its highest-degree term is $P_dP_d^*\ne0$ in
the free algebra. The change from $x_\ell,x_\ell^*$ to
$s_{\ell,1},s_{\ell,2}$ is invertible, and semicircular variables
satisfy no polynomial relations. Thus the nonconstant scalar case of
\cite{SS} applies. Use the convention
\[
 g^z(w)=\int(x-w)^{-1}\nu^z(\dd x),
\]
whose imaginary part is positive in the upper half-plane. The
non-atomicity and the local Puiseux expansion imply that, for some
$\theta>0$,
\begin{equation}\label{eq:free-small}
 \nu^z([0,x])\le Cx^\theta,\qquad 0<x\le x_*.
\end{equation}
Indeed, the Cauchy bound $|g^z(i\eta)|\le\eta^{-1}$
excludes Puiseux powers below $-1$, and the absence of an atom at zero
excludes a nonzero $w^{-1}$ term. After decreasing $\theta$ to be
less than $1$, it follows that
$\Im g^z(i\eta)\le C\eta^{\theta-1}$. Since
\[
 \nu^z([0,x])\le2x\Im g^z(ix),
\]
we obtain \eqref{eq:free-small}. In particular the negative logarithm
is integrable against $\nu^z$.

Let $G=(G_1,\ldots,G_n)$ be a family of independent normalized complex
Ginibre matrices. Apply the GUE resolvent estimate to $Q_z$, using their
representation by pairs of independent GUE matrices. Choose
$0<a_0<b_0$ sufficiently small that
$N^{-2}\eta^{-b_1}\le C\eta^{\theta-1}$ for
$\eta\ge N^{-a_0}$. The preceding mass inequality gives
\begin{equation}\label{eq:gaussian-count}
 \EE N^{-1}\#\{j:s_j(P(G,G^*)-zI)^2\le x\}
 \le Cx^\theta,\qquad N^{-a_0}\le x\le1.
\end{equation}
We have enlarged $C$ to include $x$ away from zero.

To compare entry laws, we use the affine pencil: one atom changes it
by bounded rank, which controls resolvent derivatives. We first transfer
the Gaussian count to this pencil. Let $h$ be its coefficient size and
define its Hermitization and expected Cauchy transform by
\[
 \mathcal H_z(X)=
 \begin{pmatrix}0&\mathcal L_z(X)\\\mathcal L_z(X)^*&0\end{pmatrix},
 \qquad
 m_{X,N}(w)=
 \EE\tr_{2hN}(\mathcal H_z(X)-w)^{-1}.
\]
The Gaussian norm estimate follows from a net argument. A $1/4$-net of the unit sphere of $\C^N$ has at most
$9^{2N}$ points, and nets for both vectors give
\[
 \norm{G_\ell}\le2\max_{u,v\text{ in the nets}}|u^*G_\ell v|.
\]
For fixed unit vectors, $u^*G_\ell v$ is proper complex Gaussian of
variance $1/N$, with tail $e^{-Nt^2}$. A union bound over at most
$9^{4N}$ pairs and the fixed number of input matrices shows that
\[
 \PP\{\max_\ell\norm{G_\ell}>C_0\}\le Ce^{-cN}
\]
for sufficiently large fixed $C_0$. On the complementary event,
$\mathcal V_z,\mathcal W_z$ and their inverses have fixed norm bounds. Set
$\Delta_z=\diag(P(G,G^*)-zI,I)$. The identity
$\Delta_z=\mathcal V_z\mathcal L_z\mathcal W_z$ and min--max give
$s_j(\Delta_z)\le C_1s_j(\mathcal L_z)$. For $C_1t<1$, the identity
block of $\Delta_z$ contributes no singular values at most $C_1t$,
and hence
\[
 \#\{s_j(\mathcal L_z)\le t\}
 \le\#\{s_j(P(G,G^*)-zI)\le C_1t\}.
\]
Hermitization counts every singular value twice and normalizes by
$2hN$. Thus \eqref{eq:gaussian-count} bounds its expected mass in
$[-t,t]$ by $Ct^{2\theta}+Ce^{-cN}$ whenever
$(C_1t)^2\ge N^{-a_0}$ and $t$ is small. We may therefore choose
$0<\kappa<\min(2\theta,1)$ and $0<a_1<a_0/2$ so that this expected
mass is at most $Ct^\kappa$ for $N^{-a_1}\le t\le1$.

For $\eta\ge N^{-a_1}$, integrate the kernel
$\eta/(\lambda^2+\eta^2)$ first on $|\lambda|\le\eta$, and then on
the shells $2^j\eta<|\lambda|\le2^{j+1}\eta\le1$. The first interval
contributes at most $C\eta^{\kappa-1}$. On the $j$th shell the
contribution is at most
\[
 \frac{\eta}{(2^j\eta)^2}\,C(2^{j+1}\eta)^\kappa
 \le C\eta^{\kappa-1}2^{j(\kappa-2)}.
\]
The region bounded away from zero contributes $O(\eta)$; a possible
last shell meeting $1$ changes only the constant. Summing gives
\begin{equation}\label{eq:gaussian-linear}
 \Im m_{G,N}(i\eta)\le C\eta^{\kappa-1},
 \qquad N^{-a_1}\le\eta\le1.
\end{equation}
The counting bounds apply throughout the summation, since all scales
are at least $\eta$.

\subsection{Resolvent comparison for general entry laws}

We transfer the Gaussian pencil bound in two steps. Replacing each atom
by a Gaussian with the same full real covariance uses the fractional
moment assumption; Gaussian interpolation then removes the
pseudovariance. Both errors must be small enough at the counting scale.
Write the Hermitization as
\begin{equation}\label{eq:Hermitian-atoms}
 \mathcal H_z(X)=D_z\otimes I+
 \frac1{\sqrt N}\sum_{\ell,a,b}
  \left(\xi_{\ell,ab}\mathcal C_\ell\otimes E_{ab}
       +\overline{\xi_{\ell,ab}}\mathcal C_\ell^*\otimes E_{ba}\right),
 \qquad
 \mathcal C_\ell=\begin{pmatrix}0&U_\ell\\V_\ell^*&0\end{pmatrix},
\end{equation}
where $D_z=\left(\begin{smallmatrix}0&K_z\\K_z^*&0\end{smallmatrix}\right)$.
Indeed, the upper-right blocks contributed by one atom are
$\xi U_\ell E_{ab}$ and $\bar\xi V_\ell E_{ba}$, the two terms of
$\mathcal L_z$; the lower-left blocks are their adjoints.

Let $\gamma_\rho$ be a centered complex Gaussian with real covariance
\eqref{eq:Sigma}, and let $G_\rho$ have independent entries
$\gamma_\rho/\sqrt N$. Thus $\EE|\gamma_\rho|^2=1$ and
$\EE\gamma_\rho^2=\rho$. The choice $\rho=1$ is a standard real
Gaussian, and $G=G_0$ is the standard complex Ginibre family. We compare
$X$ with $G_\rho$, then $G_\rho$ with $G$.

\begin{lemma}[Fractional comparison of whole atoms]\label{lem:fractional}
For $0<\alpha\le1$, if $\EE|\xi|^{2+\alpha}<\infty$ and
$\Im w>0$, then, in either input class,
\[
 |m_{X,N}(w)-m_{G_\rho,N}(w)|
 \le CN^{-\alpha/2}(\Im w)^{-3-\alpha}.
\]
No density assumption is needed for this lemma.
\end{lemma}
\begin{proof}
Replace the $nN^2$ independent entries one at a time, treating the
real and imaginary parts of each entry as a single atom. At each replacement, condition on the other entries. The Hermitian
background $H^{(0)}$ is then independent of the selected atom. For
$v=(v_1,v_2)$ and $\chi=v_1+iv_2$, write its perturbation as
\[
 A(v)=N^{-1/2}\bigl(\chi\mathcal C_\ell\otimes E_{ab}
                 +\bar \chi\mathcal C_\ell^*\otimes E_{ba}\bigr).
\]
It is Hermitian, has rank at most $4h$, and satisfies
$\norm{A(v)}\le C\norm v/\sqrt N$. For real atoms take $v_2=0$.
Set $f_v(t)=\tr_{2hN}(H^{(0)}+tA(v)-w)^{-1}$ and $\eta=\Im w$.
With $\mathcal R_t=(H^{(0)}+tA(v)-w)^{-1}$, differentiation gives
\[
 f_v^{(j)}(t)=(-1)^jj!\tr_{2hN}[\mathcal R_t(A(v)\mathcal R_t)^j].
\]
For $j\ge1$, the product has bounded rank and
$\norm{\mathcal R_t}\le\eta^{-1}$. The normalized trace contributes a factor
$N^{-1}$, so
\[
 |f_v''(t)|\le CN^{-2}\eta^{-3}\norm v^2,\qquad
 |f_v^{(3)}(t)|\le CN^{-5/2}\eta^{-4}\norm v^3.
\]

The second-order Taylor remainder in the two real coordinates is
\[
 \mathrm{Rem}_2(v)=f_v(1)-f_v(0)-f_v'(0)-\tfrac12 f_v''(0)
       =\int_0^1(1-t)(f_v''(t)-f_v''(0))\dd t.
\]
Bounding the integrand by the second derivatives, or integrating the
third derivative, gives
\begin{align*}
 |\mathrm{Rem}_2(v)|
 &\le C\min\{N^{-2}\eta^{-3}\norm v^2,
                 N^{-5/2}\eta^{-4}\norm v^3\}\\
 &\le CN^{-2-\alpha/2}\eta^{-3-\alpha}\norm v^{2+\alpha},
\end{align*}
where the last step uses $\min(a,b)\le a^{1-\alpha}b^\alpha$.
The constant term does not depend on $v$; the next two terms are
respectively linear in $v_1,v_2$ and quadratic in these coordinates. Their expectations agree for $\xi$
and $\gamma_\rho$, because the means and the entire real covariance
agree, including the mixed $v_1v_2$ entry. This remains true for the
rank-one covariance of class (A). Summing the expected remainders over
$nN^2$ replacements proves the claim. Every intermediate background
is Hermitian, so the resolvent bound holds regardless of the other entries.
\end{proof}

The following covariance comparison develops the argument of
\cite[Lemma~6.11]{Han}.

\begin{lemma}[Removing arbitrary Gaussian pseudovariance]\label{lem:rho}
For $|\rho|\le1$ and $\Im w>0$,
\[
 |m_{G_\rho,N}(w)-m_{G,N}(w)|
 \le CN^{-1}(\Im w)^{-3}.
\]
The constant may be chosen uniformly over $|\rho|\le1$.
\end{lemma}
\begin{proof}
Let $Y_\rho,Y_0$ denote the centered Gaussian parts of
\eqref{eq:Hermitian-atoms}. For deterministic $M$, define the fixed
linear covariance maps
$\mathcal K_\rho(M)=\EE[Y_\rho M Y_\rho]$ and
$\mathcal K_0(M)=\EE[Y_0MY_0]$. We apply these maps pointwise when
$M$ is random. Terms involving distinct entries vanish by independence
and centering.
The two terms pairing an entry with its conjugate agree because the
absolute variances are one. The other two terms give
\begin{align}\label{eq:covariance-difference}
 (\mathcal K_\rho-\mathcal K_0)(M)
 =\frac1N\sum_{\ell,a,b}\bigl[
 &\rho(\mathcal C_\ell\otimes E_{ab})M
                       (\mathcal C_\ell\otimes E_{ab})\notag\\
 &+\bar\rho(\mathcal C_\ell^*\otimes E_{ba})M
                       (\mathcal C_\ell^*\otimes E_{ba})\bigr].
\end{align}
For an $N\times N$ matrix $Q$,
$\sum_{a,b}E_{ab}QE_{ab}=Q^{\mathsf T}$. Thus the sums in
\eqref{eq:covariance-difference} transpose each spatial block of
$M$, with fixed coefficient matrices on either side. Every such
block has norm at most $\norm M$, transpose preserves its norm, and
a $2h\times2h$ block array has norm at most $2h$ times its largest
block norm. It follows that
\begin{equation}\label{eq:covariance-map-bound}
 \norm{(\mathcal K_\rho-\mathcal K_0)(M)}
 \le CN^{-1}\norm M,
\end{equation}
uniformly in $|\rho|\le1$.

Interpolate independent copies by
$H_t=D_z\otimes I+\sqrt tY_\rho+\sqrt{1-t}Y_0$, and write
$\mathcal R_t=(H_t-w)^{-1}$. Gaussian covariance differentiation gives
\[
 \frac{\dd}{\dd t}\EE \mathcal R_t
 =\EE\bigl[\mathcal R_t(\mathcal K_\rho-\mathcal K_0)(\mathcal R_t)\mathcal R_t\bigr].
\]
To verify the identity, choose a real Hermitian coefficient basis $H_r$
and let $\Delta_{rs}$ be the difference of its real covariances.
The derivative of the Gaussian expectation is
$\frac12\sum_{r,s}\Delta_{rs}\EE\partial_r\partial_s\mathcal R$. Since
\[
 \partial_r\partial_s\mathcal R=\mathcal R H_r\mathcal R H_s\mathcal R+\mathcal R H_s\mathcal R H_r\mathcal R,
\]
symmetry of $\Delta$ combines the two terms into the covariance map
in the preceding display. By
\eqref{eq:covariance-map-bound} and $\norm{\mathcal R_t}\le\eta^{-1}$, the
derivative has norm at most $CN^{-1}\eta^{-3}$. Integrating in $t$
and taking normalized traces proves the lemma. Dominated convergence
using the same resolvent bound gives continuity at $t=0,1$. The
argument uses no inverse covariance, so it includes the singular
Gaussian endpoint $\rho=1$.
\end{proof}

\subsection{Transferring the count from the pencil to the polynomial}\label{subsec:mesoscopic-count}

Choose a polynomial scale at which the comparison errors are dominated
by the Gaussian bound. Set $\alpha=\min(\varepsilon,1)$ and let
$\widehat\mu_{X,N}$ be the eigenvalue law of $\mathcal H_z(X)$.
The two comparison lemmas and
\eqref{eq:gaussian-linear}, together with
$\EE\widehat\mu_{X,N}([-\eta,\eta])\le2\eta\Im m_{X,N}(i\eta)$,
give
\[
 \EE\widehat\mu_{X,N}([-\eta,\eta])
 \le C\bigl(\eta^\kappa+
 N^{-\alpha/2}\eta^{-2-\alpha}+N^{-1}\eta^{-2}\bigr).
\]
Choose $a_2>0$ so that
\[
 a_2<a_1,\qquad a_2(\kappa+2+\alpha)<\alpha/2,\qquad
 a_2(\kappa+2)<1.
\]
For $N^{-a_2}\le\eta\le1$, each error is then at most
$C\eta^\kappa$.

To reach $P-zI$, the Gaussian argument used bounded elimination factors.
Under the present moment assumptions, truncate these factors and bound
the discarded rank in expectation, without an operator-norm assumption
on $X$. Write
\[
 \mathcal L_z=A_{\mathrm L}\Delta_zA_{\mathrm R},\qquad
 A_{\mathrm L}=\mathcal V_z^{-1},\quad A_{\mathrm R}=\mathcal W_z^{-1},\quad
 \Delta_z=\diag(P-zI,I).
\]
The matrices $A_{\mathrm L},A_{\mathrm R}$ are fixed matrix-valued $*$-polynomials. Set their singular
values above $M$ to zero, obtaining $A_{\mathrm L,M},A_{\mathrm R,M}$ of norm at most $M$.
If $r_{\mathrm L},r_{\mathrm R}$ are the respective rank losses, then
\[
 \rank(\mathcal L_z-A_{\mathrm L,M}\Delta_zA_{\mathrm R,M})\le r_{\mathrm L}+r_{\mathrm R},\qquad
 s_j(A_{\mathrm L,M}\Delta_zA_{\mathrm R,M})\le M^2s_j(\Delta_z).
\]
For $0<x<1$, the identity block contributes nothing to the count of
$s_j(\Delta_z)\le\sqrt x$. Every singular value counted there is
therefore counted in $A_{\mathrm L,M}\Delta_zA_{\mathrm R,M}$ at threshold $M^2\sqrt x$.
Rank interlacing with $\mathcal L_z$ now gives, for
$\mathcal N_N^z(x)=N^{-1}\#\{j:s_j(P-zI)^2\le x\}$,
\begin{equation}\label{eq:rank-transfer}
 \mathcal N_N^z(x)\le
 h\widehat\mu_{X,N}([-M^2\sqrt x,M^2\sqrt x])
       +(r_{\mathrm L}+r_{\mathrm R})/N,\qquad 0<x<1.
\end{equation}
The factor $h$ follows because Hermitization counts each singular
value twice and uses normalization $2hN$.

By Lemma~\ref{sstyle:lem:rank-cap},
$\EE(r_{\mathrm L}+r_{\mathrm R})\le CNM^{-r_*}$ for some fixed $r_*>0$. Choose
$M=x^{-1/8}$. Then $M^2\sqrt x=x^{1/4}$, so the preceding pencil
count and \eqref{eq:rank-transfer} yield
\begin{equation}\label{eq:mesoscopic}
 \EE \mathcal N_N^z(x)\le Cx^\vartheta,\qquad
 N^{-a}\le x\le x_*,\qquad
 \vartheta=\min(\kappa/4,r_*/8)>0,
\end{equation}
where any sufficiently small $a<4a_2$ is admissible. The count is
unconditional, since the rank errors have been estimated in expectation.

\subsection{Logarithmic tails and pointwise convergence of potentials}\label{subsec:log-potentials}

For assertion (ii), combine the small-value count with the least
singular value cutoff to control the negative logarithmic tail. The
positive tail needs only second moments and the rank bound. Let $d=\deg P\ge1$. Lemma~\ref{sstyle:lem:rank-cap} at threshold
$\sqrt x$ gives $\EE\nu_N^z((x,\infty))\le Cx^{-1/d}$. Integrating
this estimate yields
\begin{align}\label{eq:upper-log}
 \sup_N\EE\int_{(T,\infty)}\log x\,\nu_N^z(\dd x)
 &\le C T^{-1/d}\log T+C\int_T^\infty x^{-1-1/d}\dd x\notag\\
 &\le C T^{-1/d}(\log T+d),\qquad T\ge1.
\end{align}

The count alone does not exclude a few extremely small singular values.
For the lower tail, Corollary~\ref{cor:poly} gives a fixed $B$, sufficiently
large that
\[
 \mathcal F_N=\{\smin(P-zI)\ge N^{-B}\},\qquad
 \PP(\mathcal F_N^c)\longrightarrow0.
\]
Set $x_N=N^{-a}$ with $a$ from \eqref{eq:mesoscopic}, enlarging $B$
so that $2B>a$. On $\mathcal F_N$, every squared singular value is
at least $N^{-2B}$, and hence
\[
 \int_{(0,x_N]}|\log x|\,\nu_N^z(\dd x)
 \le2B\log N\,\mathcal N_N^z(x_N).
\]
Its expectation after multiplication by $\one_{\mathcal F_N}$ is at
most $CN^{-a\vartheta}\log N=o(1)$. On the remaining interval,
integration by parts gives, for $x_N\le\delta\le\min(x_*,1)$,
\begin{align*}
 \EE\int_{(x_N,\delta]}|\log x|\,\nu_N^z(\dd x)
 &\le |\log\delta|\,\EE \mathcal N_N^z(\delta)
            +\int_{x_N}^{\delta}\frac{\EE \mathcal N_N^z(x)}x\dd x\\
 &\le C\delta^\vartheta
                (|\log\delta|+\vartheta^{-1}).
\end{align*}
In particular, for every $t>0$,
\[
 \lim_{\delta\downarrow0}\limsup_{N\to\infty}
 \PP\left\{\int_{(0,\delta]}|\log x|\,\nu_N^z(\dd x)>t\right\}=0.
\]
Indeed, Markov's inequality on $\mathcal F_N$ gives the bound after
adding $\PP(\mathcal F_N^c)$. The error $q_N$ in the least-singular-value
estimate appears only in this probability and acquires no factor $\log N$.

Assertion (i) now gives convergence for the logarithm truncated
above and below. The bound \eqref{eq:free-small} removes its lower
truncation in the free law, whose upper tail vanishes by compact
support. The preceding lower-tail estimate and \eqref{eq:upper-log}
remove the two truncations for the random laws. Recalling that
$\frac12\int\log x\dd\nu_N^z=N^{-1}\log|\det(P-zI)|$, we obtain
\begin{equation}\label{eq:potential}
 U_N(z):=N^{-1}\log|\det(P(X,X^*)-zI)|
 \ \xrightarrow{\PP}
 U(z):=\tau(\log|P(c,c^*)-z|),\qquad z\ne P(0).
\end{equation}

\subsection{Spatial integrability and Brown measure convergence}\label{subsec:spatial}

To pass to the Laplacian, we upgrade pointwise potential convergence to
convergence after integrating the absolute difference over compact planar
sets. A uniform spatial square-integral bound gives the required uniform
integrability. We obtain it from the logarithmic kernel and eigenvalue
product inequalities, starting with the consequence of
\eqref{eq:upper-count},
\[
 \sup_N\EE\frac1N\sum_j(\log^+s_j(P))^2<\infty.
\]
To apply the eigenvalue product inequalities, write $a_j=\log^+|\lambda_j(P)|$ and $b_j=\log^+s_j(P)$ in
decreasing order. For a fixed $k$, let $r\le k$ be the number of
positive eigenvalue logarithms among the first $k$. Then
\[
 \sum_{j\le k}a_j=\log\prod_{j\le r}|\lambda_j(P)|
 \le\log\prod_{j\le r}s_j(P)
 \le\sum_{j\le k}b_j.
\]
Thus, for $t\ge0$,
\[
 \sum_j(a_j-t)_+
 =\max_{0\le k\le N}\left(\sum_{j\le k}a_j-kt\right)
 \le\sum_j(b_j-t)_+.
\]
Integrating against $2\dd t$ proves
$\sum_ja_j^2\le\sum_jb_j^2$. In particular the eigenvalue measures
are tight in probability, and their expected squared positive
logarithms are uniformly bounded.

For a compact planar set $K$, the elementary kernel estimate
\[
 \int_K|\log|z-\lambda||^2\dd^2z
 \le C_K(1+(\log^+|\lambda|)^2)
\]
follows by integrating $|\log|z-\lambda||^2$ on $|z-\lambda|\le1$ and bounding
the positive logarithm on the remaining part of $K$. Since
$U_N(z)=N^{-1}\sum_j\log|z-\lambda_j(P)|$, Jensen's inequality and
the preceding eigenvalue bound give
\begin{equation}\label{sstyle:eq:spatial-L2}
 \sup_N\EE\int_K|U_N(z)|^2\dd^2z<\infty.
\end{equation}
Realize all matrix sizes on a common product probability space. For any $t>0$,
\eqref{eq:potential} and bounded convergence in the planar variable
give
\[
 \int_K\PP\{|U_N(z)-U(z)|>t\}\dd^2z\longrightarrow0.
\]
The single excluded point $P(0)$ has zero area. Hence $U_N$ converges
to $U$ in measure on the finite product space. Take a subsequence
converging almost everywhere there. Fatou's lemma and
\eqref{sstyle:eq:spatial-L2} give $U\in L^2(K)$, and the same bound
makes $U_N-U$ uniformly integrable on the product space. It follows
that
\[
 \EE\int_K|U_N(z)-U(z)|\dd^2z\longrightarrow0.
\]
Thus the potentials converge locally in expected spatial $L^1$.

Finally $\Delta U_N/(2\pi)$ is the empirical eigenvalue measure and
$\Delta U/(2\pi)$ is the Brown measure. For a smooth compactly
supported $\phi$,
\[
 \left|\int\phi\dd\mu_{P(X,X^*)}-\int\phi\dd\mu_{P(c,c^*)}\right|
 \le\frac{\norm{\Delta\phi}_\infty}{2\pi}
       \int_{\operatorname{supp}\phi}|U_N-U|\dd^2z,
\]
which tends to zero in probability. Eigenvalue tightness extends the
convergence to bounded continuous tests. Constant $P$ is immediate,
and Theorem~\ref{thm:brown} follows.

\subsection{Spectral transfer for Wigner inputs}\label{w:sec:spectral}\label{subsec:wigner-spectral}

For Theorem~\ref{thm:wigner-brown}, the limit is semicircular. The tail
and spatial arguments apply once we compare Wigner atoms: an
off-diagonal atom occupies conjugate positions, and the diagonal
replacement must use only second moments. Off-diagonal laws have finite
$(2+\varepsilon)$ moments for some $\varepsilon>0$; diagonal laws need only
second moments. The density assumptions enter through
Theorem~\ref{thm:wigner-pencil}. The entry laws may differ among input
matrices. Define
\[
 \rho_\ell=\EE\zeta_{\ell,12}^2,\qquad
 \tau_\ell^2=\EE d_{\ell,1}^2,\qquad
 \alpha=\min(\varepsilon,1).
\]
This includes real atoms ($\rho_\ell=1$) and correlated planar atoms.

Fix a nonconstant scalar polynomial $P$ in the Hermitian letters and
$z\ne P(0)$. Identify each letter with its adjoint in
Lemma~\ref{lem:linearization} and write the resulting pencil as
\[
 \mathcal L_z(H)=K_z\otimes I_N+
                 \sum_{\ell=1}^nC_{\ell,z}\otimes H_\ell.
\]
If $h$ is its coefficient size, its Hermitization is
\begin{equation}\label{w:eq:hermitization}
 \mathcal H_z(H)=D_z\otimes I_N+
                    \sum_\ell B_{\ell,z}\otimes H_\ell,
 \quad
 D_z=\begin{pmatrix}0&K_z\\K_z^*&0\end{pmatrix},\quad
 B_{\ell,z}=\begin{pmatrix}0&C_{\ell,z}\\C_{\ell,z}^*&0\end{pmatrix}.
\end{equation}
All $B_{\ell,z}$ are Hermitian. For $\Im w>0$, define
\[
 m_{H,N}(w)=
 \EE\tr_{2hN}(\mathcal H_z(H)-w)^{-1},\qquad
 \eta=\Im w.
\]
The fixed shift $z$ is suppressed in this notation.

\paragraph{Comparison of Wigner atoms.}
Let $H^{\rm match}$ have independent centered Gaussian atoms matching
the Wigner real covariances (off-diagonal absolute variance one,
pseudovariance $\rho_\ell$, diagonal variance $\tau_\ell^2$).
Let $G^{\mathrm{GUE}}_1,\ldots,G^{\mathrm{GUE}}_n$ be independent normalized GUE matrices, with
both off-diagonal absolute variance and diagonal variance $1/N$.

\begin{lemma}[Wigner replacement and covariance comparison]
\label{w:lem:comparison}
For every $\Im w>0$,
\begin{align}
 |m_{H,N}(w)-m_{H^{\rm match},N}(w)|
 &\le C N^{-\alpha/2}\eta^{-3-\alpha}
          +C N^{-1}\eta^{-3},\label{w:eq:entry-comparison}\\
 |m_{H^{\rm match},N}(w)-m_{G^{\mathrm{GUE}},N}(w)|
 &\le C N^{-1/2}\eta^{-3}.\label{w:eq:gaussian-comparison}
\end{align}
Constants may depend on the fixed coefficients and atom laws. No
density assumption is needed for these two estimates.
\end{lemma}
\begin{proof}
At $a<b$, keep the off-diagonal atom $\chi=v_1+iv_2$ together with
its conjugate reverse entry. Its Hermitized perturbation is
\[
 A(v)=N^{-1/2}B_{\ell,z}\otimes
                    (\chi E_{ab}+\bar\chi E_{ba}).
\]
This is Hermitian, has bounded rank and norm at most
$C\norm v/\sqrt N$, and is independent of the background after conditioning
on the other atoms. Apply the Taylor calculation of
Lemma~\ref{lem:fractional} to the two real coordinates together; its remainder is
$CN^{-2-\alpha/2}\eta^{-3-\alpha}\norm v^{2+\alpha}$.
Matching the entire real covariance, including its mixed
real-imaginary entry, makes the quadratic expectations agree; centering
does the same for the linear terms. Summing over $O(N^2)$ atoms gives an error of
$CN^{-\alpha/2}\eta^{-3-\alpha}$.

For a diagonal atom the perturbation is
$N^{-1/2}dB_{\ell,z}\otimes E_{aa}$. Expand only to first order.
The first-order expectations agree by centering, and the second
resolvent derivative bounds the expected remainder by
$CN^{-2}\eta^{-3}\EE|d|^2$, together with the corresponding Gaussian
term. Summing over the $O(N)$ diagonal atoms gives
$CN^{-1}\eta^{-3}$. This proves \eqref{w:eq:entry-comparison} without
a moment above two on the diagonal.

For the Gaussian comparison, denote the centered Hermitized parts by $Y_{\rm match},Y_{\mathrm{GUE}}$, and define fixed linear maps
$\mathcal K_\bullet(M)=\EE[Y_\bullet M Y_\bullet]$ for deterministic
$M$. The off-diagonal terms pairing an entry with its conjugate cancel
because the absolute variances match. The remaining terms are
\begin{align}\label{w:eq:covariance-difference}
 (\mathcal K_{\rm match}-\mathcal K_{\mathrm{GUE}})(M)
 =\frac1N\sum_{\ell,a<b}\bigl[&
 \rho_\ell(B_{\ell,z}\otimes E_{ab})M
                         (B_{\ell,z}\otimes E_{ab})\notag\\
 &+\bar\rho_\ell(B_{\ell,z}\otimes E_{ba})M
                         (B_{\ell,z}\otimes E_{ba})\bigr]\notag\\
 +\frac1N\sum_{\ell,a}(\tau_\ell^2-1)&
 (B_{\ell,z}\otimes E_{aa})M(B_{\ell,z}\otimes E_{aa}).
\end{align}
To bound this map, consider one spatial $N\times N$ block $Q$ of
$M$. The two off-diagonal sums act as
\[
 \sum_{a<b}E_{ab}QE_{ab}
      =\operatorname{tri}_{+}(Q^{\mathsf T}),\qquad
 \sum_{a<b}E_{ba}QE_{ba}
      =\operatorname{tri}_{-}(Q^{\mathsf T}),
\]
where $\operatorname{tri}_{+}$ and $\operatorname{tri}_{-}$ denote
strict upper and lower triangular truncation. Their weights are
$\rho_\ell$ and $\bar\rho_\ell$, so for nonreal pseudovariance the
sum need not be a scalar multiple of the full transpose. We bound each part by
\[
 \norm{\operatorname{tri}_{\pm}(Q^{\mathsf T})}
 \le\norm{\operatorname{tri}_{\pm}(Q^{\mathsf T})}_{\HS}
 \le\norm Q_{\HS}\le\sqrt N\norm Q.
\]
The diagonal sum is $\diag(Q_{11},\ldots,Q_{NN})$ and has norm at
most $\norm Q$. Since the coefficient block dimensions are fixed,
\eqref{w:eq:covariance-difference} gives
\begin{equation}\label{w:eq:covariance-map-bound}
 \norm{(\mathcal K_{\rm match}-\mathcal K_{\mathrm{GUE}})(M)}
 \le CN^{-1/2}\norm M.
\end{equation}

Apply the Gaussian interpolation identity from
Lemma~\ref{lem:rho} to independent $Y_{\rm match},Y_{\mathrm{GUE}}$.
By \eqref{w:eq:covariance-map-bound}, its derivative has norm at most
$CN^{-1/2}\eta^{-3}$. Integrating and taking the normalized trace gives
\eqref{w:eq:gaussian-comparison}. The resolvent bound also gives continuity at the endpoints, including
singular real Gaussian covariances.
\end{proof}

\paragraph{The polynomial count.}
We combine the Wigner comparison with the GUE reference bound and the
expected rank loss from truncating the elimination factors.

\begin{lemma}[Wigner mesoscopic singular-value count]
\label{w:lem:count}
Let
\[
 \mathcal N_{H,N}^z(x)=N^{-1}\#\{j:s_j(P(H)-zI_N)^2\le x\}.
\]
There are fixed $a,\vartheta,x_*,C>0$, depending on $P,z$ and the
atom laws, such that
\begin{equation}\label{w:eq:mesoscopic}
 \EE \mathcal N_{H,N}^z(x)\le Cx^\vartheta,
 \qquad N^{-a}\le x\le x_*.
\end{equation}
This estimate is unconditional.
\end{lemma}
\begin{proof}
The argument of Section~\ref{subsec:gaussian-reference} applies to independent
normalized GUE matrices. The squared singular-value limit is $\nu_{\rm sc}^z$, the law of
$Q_z=(P(s)-z)(P(s)-z)^*$ for free semicircular $s_1,\ldots,s_n$.
The highest homogeneous term $P_dP_d^*\ne0$ makes $Q_z$ nonconstant,
so \cite[Theorem~1.1 and Corollary~1.2]{SS} and the Puiseux argument
for \eqref{eq:free-small} give
\begin{equation}\label{w:eq:free-small}
 \nu_{\rm sc}^z([0,x])\le Cx^\theta,
 \qquad 0<x\le x_*,\qquad \theta>0.
\end{equation}
In particular its negative logarithm is integrable. Applying
\cite[Lemma~5.5.4]{AGZ} directly to $Q_z$, choosing the cutoff as in \eqref{eq:gaussian-count}, gives
\begin{equation}\label{w:eq:gaussian-count}
 \EE N^{-1}\#\{j:s_j(P(G^{\mathrm{GUE}})-zI)^2\le x\}
 \le Cx^\theta,
 \qquad N^{-a_0}\le x\le1
\end{equation}
for some $a_0>0$. For GUE, the real and imaginary parts of
$u^*G^{\mathrm{GUE}}_\ell v$ have variance at most $C/N$ for fixed complex unit
vectors $u,v$. The net argument of
Section~\ref{subsec:gaussian-reference} gives
$\PP\{\max_\ell\norm{G^{\mathrm{GUE}}_\ell}>C_0\}\le Ce^{-cN}$.
On the complementary event the elimination factors are bounded, so
the min--max argument and shell summation in \eqref{eq:gaussian-linear} give
\begin{equation}\label{w:eq:gaussian-linear}
 \Im m_{G^{\mathrm{GUE}},N}(i\eta)\le C\eta^{\kappa-1},
 \qquad N^{-a_1}\le\eta\le1,
\end{equation}
where $0<\kappa<\min(2\theta,1)$ and $a_1>0$.

Let $\widehat\mu_{H,N}$ be the eigenvalue law of $\mathcal H_z(H)$.
Lemma~\ref{w:lem:comparison}, \eqref{w:eq:gaussian-linear}, and the
resolvent mass inequality give
\[
 \EE\widehat\mu_{H,N}([-\eta,\eta])
 \le2\eta\Im m_{H,N}(i\eta)
 \le C\{\eta^\kappa+
        N^{-\alpha/2}\eta^{-2-\alpha}
        +N^{-1/2}\eta^{-2}\}.
\]
The last term includes the smaller error from the diagonal atoms. Choose $a_2>0$ satisfying
\begin{equation}\label{w:eq:scale-choice}
 a_2<a_1,\qquad
 a_2(\kappa+2+\alpha)<\alpha/2,\qquad
 a_2(\kappa+2)<1/2.
\end{equation}
Then the last bound is at most $C\eta^\kappa$ for
$N^{-a_2}\le\eta\le1$.

The Wigner second moments satisfy $\EE\norm{H_\ell}_{\HS}^2=N-1+\tau_\ell^2\le CN$.
Lemma~\ref{sstyle:lem:rank-cap} thus applies to every fixed matrix-valued
polynomial $Q$ in the Hermitian letters, of degree at most $D\ge1$, and gives
\begin{equation}\label{w:eq:upper-count}
 \EE\#\{j:s_j(Q(H))>M\}\le C_QNM^{-2/D},\qquad M\ge1.
\end{equation}
Apply \eqref{eq:rank-transfer} and the argument of
Section~\ref{subsec:mesoscopic-count} to the polynomial factors $A_{\mathrm L}=\mathcal V_z^{-1}$ and $A_{\mathrm R}=\mathcal W_z^{-1}$: deleting their
singular values above $M$ has expected rank loss at most $CNM^{-r_*}$
for some $r_*>0$. Taking $M=x^{-1/8}$ proves
\eqref{w:eq:mesoscopic} with
$\vartheta=\min(\kappa/4,r_*/8)>0$ and sufficiently small $a<4a_2$.
The expected rank bounds make the count unconditional.
\end{proof}

\paragraph{The logarithmic limit.}
To complete the logarithmic integral, combine this count with weak
convergence and the Wigner least singular value cutoff.

\begin{lemma}[Finite-variance weak singular convergence for Wigner inputs]
\label{w:lem:weak}
Under only the finite second moments of the Wigner model, the squared
singular-value law of $P(H)-zI$ converges weakly in probability to
$\nu_{\rm sc}^z$.
\end{lemma}
\begin{proof}
For bounded laws, use the closed-walk expansion of
Section~\ref{subsec:weak-log}, now indexed by unordered edges labeled by $\ell$.
A surviving length-$k$ pattern has $v\le e+1\le1+k/2$.
Equality forces a tree traversed once in each direction, with each
edge contributing $\EE|\zeta_{\ell,12}|^2=1$ regardless of
$\rho_\ell$. Loops (diagonal atoms) and parallel edges with different matrix labels lose
a free index. Leaf removal gives the semicircular noncrossing
pairings with matching matrix labels. The covariance count for connected walks is $O(N^{-1})$;
compact support and moment determinacy give weak convergence
in probability, as in Section~\ref{subsec:weak-log}.

For finite-variance laws, truncate and center the independent atoms
at $L$, normalizing only the off-diagonal absolute variances to one.
Reflecting the upper-triangular entries by conjugation gives a Hermitian
coupling $H_{\ell,L}$. With generic atoms $\zeta_\ell,d_\ell$ and their
truncated approximants on the right, it satisfies
\begin{align*}
 \EE N^{-1}\norm{H_\ell-H_{\ell,L}}_{\HS}^2
 ={}&(1-N^{-1})\EE|\zeta_\ell-\zeta_{\ell,L}|^2\\
 &+N^{-1}\EE|d_\ell-d_{\ell,L}|^2.
\end{align*}
The right side is at most $\delta_L^2\to0$, uniformly in $N,\ell$.
Diagonal variances do not affect the leading trees. Since truncation
preserves Hermitian symmetry, we may apply \eqref{sstyle:eq:truncation-test} from
Section~\ref{subsec:weak-log}. Choose $T$, then $L$, then let
$N\to\infty$; tightness from \eqref{w:eq:upper-count} proves the lemma.
\end{proof}

\begin{proof}[Proof of Theorem~\ref{thm:wigner-brown}]
Lemma~\ref{w:lem:weak} gives weak singular-value convergence. Fix $z\ne P(0)$. The pencil
$(K_z^{-1}\otimes I)\mathcal L_z(H)$ has identity constant term.
Theorem~\ref{thm:wigner-pencil} and the inverse-block identity in
Lemma~\ref{lem:linearization} give a fixed $B_z<\infty$ such that
\[
 \mathcal F_{N,z}=\{s_{\min}(P(H)-zI)\ge N^{-B_z}\},\qquad
 \PP(\mathcal F_{N,z}^c)\longrightarrow0.
\]
This cutoff uses only second moments, including on the diagonal.

Use the lower-tail argument of Section~\ref{subsec:log-potentials}
with \eqref{w:eq:mesoscopic}. Let $x_N=N^{-a}$ and enlarge $B_z$ so that $2B_z>a$.
For the squared singular-value law $\nu_{H,N}^z$ it gives
\[
 \EE\one_{\mathcal F_{N,z}}
       \int_{(0,x_N]}|\log x|\,\nu_{H,N}^z(\dd x)
 \le CN^{-a\vartheta}\log N=o(1),
\]
and, for $x_N\le\delta\le\min(x_*,1)$, the integral over
$(x_N,\delta]$ has expectation at most
$C\delta^\vartheta(|\log\delta|+\vartheta^{-1})$.
Markov's inequality on $\mathcal F_{N,z}$ bounds the lower tail after
adding the exceptional probability, without a factor $\log N$. The upper tail follows from
\eqref{w:eq:upper-count} by \eqref{eq:upper-log}. Together with
Lemma~\ref{w:lem:weak} and \eqref{w:eq:free-small}, these estimates
remove both logarithmic truncations and give
\[
 U_{H,N}(z):=N^{-1}\log|\det(P(H)-zI)|
 \xrightarrow{\PP}
 U_{\rm sc}(z):=\tau(\log|P(s)-z|),\qquad z\ne P(0).
\]

For spatial integrability, use \eqref{w:eq:upper-count} to bound the
expected squared positive logarithms of singular values, as in
Section~\ref{subsec:spatial}. The eigenvalue product inequalities,
planar kernel estimate, and uniform integrability give, for every compact $K$,
\[
 \EE\int_K|U_{H,N}(z)-U_{\rm sc}(z)|\dd^2z\longrightarrow0.
\]
Taking distributional Laplacians and using eigenvalue tightness gives
weak convergence in probability to the Brown measure. The case of constant
$P$ is immediate, completing the proof.
\end{proof}

\printbibliography[heading=bibintoc]

@book{AGZ,
  author = {Anderson, Greg W. and Guionnet, Alice and Zeitouni, Ofer},
  title = {An Introduction to Random Matrices},
  series = {Cambridge Studies in Advanced Mathematics},
  volume = {118},
  publisher = {Cambridge University Press},
  location = {Cambridge},
  year = {2010},
  doi = {10.1017/CBO9780511801334}
}

@article{AEKN,
  author = {Alt, Johannes and Erd{\H{o}}s, L{\'a}szl{\'o} and Kr{\"u}ger, Torben and Nemish, Yuriy},
  title = {Location of the spectrum of {Kronecker} random matrices},
  journaltitle = {Annales de l'Institut Henri Poincar{\'e}, Probabilit{\'e}s et Statistiques},
  volume = {55},
  number = {2},
  year = {2019},
  pages = {661--696},
  doi = {10.1214/18-AIHP894},
  eprint = {1706.08343},
  eprinttype = {arxiv}
}

@article{BSS,
  author = {Belinschi, Serban T. and {\'S}niady, Piotr and Speicher, Roland},
  title = {Eigenvalues of non-{Hermitian} random matrices and {Brown} measure of non-normal operators: {Hermitian} reduction and linearization method},
  journaltitle = {Linear Algebra and its Applications},
  volume = {537},
  year = {2018},
  pages = {48--83},
  doi = {10.1016/j.laa.2017.09.024},
  eprint = {1506.02017},
  eprinttype = {arxiv}
}

@misc{Collins,
  author = {Collins, Benoit and Kitagawa, Yu and Miyagawa, Akihiro and Yuan, Wangjun and Zhong, Ping},
  title = {{Brown} measures and spectral convergence for arbitrary {GUE} polynomials},
  year = {2026},
  eprint = {2611.12785},
  eprinttype = {arxiv}
}

@article{CGH,
  author = {Cook, Nicholas A. and Guionnet, Alice and Husson, Jonathan},
  title = {Spectrum and pseudospectrum for quadratic polynomials in {Ginibre} matrices},
  journaltitle = {Annales de l'Institut Henri Poincar{\'e}, Probabilit{\'e}s et Statistiques},
  volume = {58},
  number = {4},
  year = {2022},
  pages = {2284--2320},
  doi = {10.1214/21-AIHP1225},
  eprint = {2008.08833},
  eprinttype = {arxiv}
}

@misc{Han,
  author = {Han, Yi},
  title = {{Brown} measure convergence for the spectrum of polynomials in {Ginibre} matrices},
  year = {2026},
  eprint = {2606.01664},
  eprinttype = {arxiv}
}

@article{SS,
  author = {Shlyakhtenko, Dimitri and Skoufranis, Paul},
  title = {Freely independent random variables with non-atomic distributions},
  journaltitle = {Transactions of the American Mathematical Society},
  volume = {367},
  number = {9},
  year = {2015},
  pages = {6267--6291},
  doi = {10.1090/S0002-9947-2015-06434-4},
  eprint = {1305.1920},
  eprinttype = {arxiv}
}

@article{FK,
  author = {Fuglede, Bent and Kadison, Richard V.},
  title = {Determinant theory in finite factors},
  journaltitle = {Annals of Mathematics},
  series = {Second Series},
  volume = {55},
  number = {3},
  date = {1952},
  pages = {520--530},
  doi = {10.2307/1969645}
}

@incollection{Brown,
  author = {Brown, Lawrence G.},
  title = {Lidskii's theorem in the type {II} case},
  booktitle = {Geometric Methods in Operator Algebras (Kyoto, 1983)},
  editor = {Araki, Huzihiro and Effros, Edward G.},
  series = {Pitman Research Notes in Mathematics Series},
  number = {123},
  publisher = {Longman Scientific \& Technical},
  location = {Harlow},
  date = {1986},
  pages = {1--35}
}

@article{Girko,
  author = {Girko, Vyacheslav L.},
  title = {Circular law},
  journaltitle = {Theory of Probability \& Its Applications},
  volume = {29},
  number = {4},
  date = {1985},
  pages = {694--706},
  doi = {10.1137/1129095}
}

@article{Voiculescu,
  author = {Voiculescu, Dan},
  title = {Limit laws for random matrices and free products},
  journaltitle = {Inventiones Mathematicae},
  volume = {104},
  number = {1},
  date = {1991},
  pages = {201--220},
  doi = {10.1007/BF01245072}
}

@article{HT,
  author = {Haagerup, Uffe and Thorbj{\o}rnsen, Steen},
  title = {A new application of random matrices: {\(\operatorname{Ext}(C^*_{\mathrm{red}}(F_2))\)} is not a group},
  journaltitle = {Annals of Mathematics},
  series = {Second Series},
  volume = {162},
  number = {2},
  date = {2005},
  pages = {711--775},
  doi = {10.4007/annals.2005.162.711}
}

@article{Male,
  author = {Male, Camille},
  title = {The norm of polynomials in large random and deterministic matrices},
  journaltitle = {Probability Theory and Related Fields},
  volume = {154},
  date = {2012},
  pages = {477--532},
  doi = {10.1007/s00440-011-0375-2},
  note = {With an appendix by Dimitri Shlyakhtenko}
}

@article{Anderson,
  author = {Anderson, Greg W.},
  title = {Convergence of the largest singular value of a polynomial in independent {Wigner} matrices},
  journaltitle = {The Annals of Probability},
  volume = {41},
  number = {3B},
  date = {2013},
  pages = {2103--2181},
  doi = {10.1214/11-AOP739}
}

@article{BMS,
  author = {Belinschi, Serban T. and Mai, Tobias and Speicher, Roland},
  title = {Analytic subordination theory of operator-valued free additive convolution and the solution of a general random matrix problem},
  journaltitle = {Journal f{\"u}r die reine und angewandte Mathematik},
  volume = {732},
  date = {2017},
  pages = {21--53},
  doi = {10.1515/crelle-2014-0138}
}

@article{HMS,
  author = {Helton, J. William and Mai, Tobias and Speicher, Roland},
  title = {Applications of realizations (aka linearizations) to free probability},
  journaltitle = {Journal of Functional Analysis},
  volume = {274},
  number = {1},
  date = {2018},
  pages = {1--79},
  doi = {10.1016/j.jfa.2017.10.003}
}

@article{BBC,
  author = {Belinschi, Serban T. and Bercovici, Hari and Capitaine, Mireille},
  title = {On the outlying eigenvalues of a polynomial in large independent random matrices},
  journaltitle = {International Mathematics Research Notices},
  volume = {2021},
  number = {4},
  date = {2021},
  pages = {2588--2641},
  doi = {10.1093/imrn/rnz080}
}

@article{BC,
  author = {Bordenave, Charles and Chafa{\"i}, Djalil},
  title = {Around the circular law},
  journaltitle = {Probability Surveys},
  volume = {9},
  year = {2012},
  pages = {1--89},
  doi = {10.1214/11-PS183}
}

@article{FV,
  author = {Farrell, Brendan and Vershynin, Roman},
  title = {Smoothed analysis of symmetric random matrices with continuous distributions},
  journaltitle = {Proceedings of the American Mathematical Society},
  volume = {144},
  number = {5},
  year = {2016},
  pages = {2257--2261},
  doi = {10.1090/proc/12844}
}

@article{RVLO,
  author = {Rudelson, Mark and Vershynin, Roman},
  title = {The {Littlewood--Offord} problem and invertibility of random matrices},
  journaltitle = {Advances in Mathematics},
  volume = {218},
  number = {2},
  year = {2008},
  pages = {600--633},
  doi = {10.1016/j.aim.2008.01.010}
}

@article{RVRect,
  author = {Rudelson, Mark and Vershynin, Roman},
  title = {Smallest singular value of a random rectangular matrix},
  journaltitle = {Communications on Pure and Applied Mathematics},
  volume = {62},
  number = {12},
  year = {2009},
  pages = {1707--1739},
  doi = {10.1002/cpa.20294}
}

@article{TVK,
  author = {Tao, Terence and Vu, Van},
  title = {Random matrices: Universality of {ESDs} and the circular law},
  journaltitle = {The Annals of Probability},
  volume = {38},
  number = {5},
  year = {2010},
  pages = {2023--2065},
  note = {With an appendix by Manjunath Krishnapur},
  doi = {10.1214/10-AOP534}
}

@article{VSym,
  author = {Vershynin, Roman},
  title = {Invertibility of symmetric random matrices},
  journaltitle = {Random Structures \& Algorithms},
  volume = {44},
  number = {2},
  year = {2014},
  pages = {135--182},
  doi = {10.1002/rsa.20429}
}

@book{VDN,
  author = {Voiculescu, Dan-Virgil and Dykema, Kenneth J. and Nica, Alexandru},
  title = {Free Random Variables},
  series = {CRM Monograph Series},
  volume = {1},
  publisher = {American Mathematical Society},
  location = {Providence, RI},
  year = {1992},
  doi = {10.1090/crmm/001}
}

@book{NS,
  author = {Nica, Alexandru and Speicher, Roland},
  title = {Lectures on the Combinatorics of Free Probability},
  series = {London Mathematical Society Lecture Note Series},
  volume = {335},
  publisher = {Cambridge University Press},
  location = {Cambridge},
  year = {2006},
  doi = {10.1017/CBO9780511735127}
}

\end{document}